\documentclass[reqno,11pt]{amsart}
\usepackage{amsmath, latexsym, amsfonts, amssymb, amsthm, amscd}
\usepackage{graphics,epsf,psfrag}
\numberwithin{equation}{section}

\newtheorem{theorem}{Theorem}[section]
\newtheorem{lemma}[theorem]{Lemma}
\newtheorem{proposition}[theorem]{Proposition}
\newtheorem{cor}[theorem]{Corollary}
\newtheorem{rem}[theorem]{Remark}

\newcommand{\ind}{\mathbf{1}}
\newcommand{\E}{\mathbb{E}}
\newcommand{\R}{\mathbb{R}}
\newcommand{\Z}{\mathbb{Z}}
\newcommand{\N}{\mathbb{N}}
\renewcommand{\tilde}{\widetilde}
\renewcommand{\hat}{\widehat}

\newcommand{\cN}{{\ensuremath{\mathcal N}} }
\newcommand{\cL}{{\ensuremath{\mathcal L}} }

\newcommand{\cD}{{\ensuremath{\mathcal D}} }

\newcommand{\bP}{{\ensuremath{\mathbf P}} }
\newcommand{\bE}{{\ensuremath{\mathbf E}} }

\DeclareMathSymbol{\leqslant}{\mathalpha}{AMSa}{"36} 
\DeclareMathSymbol{\geqslant}{\mathalpha}{AMSa}{"3E} 
\DeclareMathSymbol{\eset}{\mathalpha}{AMSb}{"3F}     
\newcommand{\dd}{\,\text{\rm d}}             
\newcommand{\med}[1]{\langle #1\rangle}
\newcommand{\sumtwo}[2]{\sum_{\substack{#1 \\ #2}}} 

\newcommand{\bbE}{{\ensuremath{\mathbb E}} }

\newcommand{\bbL}{{\ensuremath{\mathbb L}} }

\newcommand{\bbP}{{\ensuremath{\mathbb P}} }

\newcommand{\ga}{\alpha}
\newcommand{\gb}{\beta}
\newcommand{\gd}{\delta}
\newcommand{\gep}{\varepsilon}       

\newcommand{\gD}{\Delta}

\newcommand{\go}{\omega}

\newcommand{\gl}{\lambda}

\makeatletter
\def\captionfont@{\footnotesize}
\def\captionheadfont@{\scshape}

\long\def\@makecaption#1#2{%
  \vspace{2mm}
  \setbox\@tempboxa\vbox{\color@setgroup
    \advance\hsize-6pc\noindent
    \captionfont@\captionheadfont@#1\@xp\@ifnotempty\@xp
        {\@cdr#2\@nil}{.\captionfont@\upshape\enspace#2}%
    \unskip\kern-6pc\par
    \global\setbox\@ne\lastbox\color@endgroup}%
  \ifhbox\@ne 
    \setbox\@ne\hbox{\unhbox\@ne\unskip\unskip\unpenalty\unkern}%
  \fi
  \ifdim\wd\@tempboxa=\z@ 
    \setbox\@ne\hbox to\columnwidth{\hss\kern-6pc\box\@ne\hss}%
  \else 
    \setbox\@ne\vbox{\unvbox\@tempboxa\parskip\z@skip
        \noindent\unhbox\@ne\advance\hsize-6pc\par}%
\fi
  \ifnum\@tempcnta<64 
    \addvspace\abovecaptionskip
    \moveright 3pc\box\@ne
  \else 
    \moveright 3pc\box\@ne
    \nobreak
    \vskip\belowcaptionskip
  \fi
\relax
}
\makeatother
\def\writefig#1 #2 #3 {\rlap{\kern #1 truecm
\raise #2 truecm \hbox{#3}}}

\newcommand{\tf}{\textsc{f}}
\newcommand{\M}{\textsc{M}}

\newcommand{\Rav}[1]{\langle R _{\substack{#1}}\rangle}

\begin{document}
\title[Hierarchical pinning, quadratic maps and disorder]{Hierarchical pinning models, quadratic maps 
\\ and quenched disorder}

\author{Giambattista Giacomin}
\address{
  Universit{\'e} Paris Diderot (Paris 7) and Laboratoire de Probabilit{\'e}s et Mod\`eles Al\'eatoires (CNRS),
U.F.R.                Math\'ematiques, Case 7012 (site Chevaleret)
             75205 Paris cedex 13, France
}
\email{giacomin\@@math.jussieu.fr}
\author{Hubert Lacoin}
\address{
  Universit{\'e} Paris Diderot (Paris 7) and Laboratoire de Probabilit{\'e}s et Mod\`eles Al\'eatoires (CNRS),
U.F.R.                Math\'ematiques, Case 7012 (site Chevaleret)
             75205 Paris cedex 13, France
}
\email{lacoin\@@math.jussieu.fr}
\author{Fabio Lucio Toninelli}
\address{CNRS and 
Laboratoire de Physique, ENS Lyon, 46 All\'ee d'Italie, 
69364 Lyon, France
}
\email{fabio-lucio.toninelli@ens-lyon.fr}
\date{\today}

\begin{abstract}
  We consider a hierarchical model of polymer pinning in presence of
  quenched disorder, introduced by 
  B. Derrida, V. Hakim and J. Vannimenus
  \cite{cf:DHV}, which can be re-interpreted as   
  an infinite dimensional dynamical system
  with random initial condition (the {\sl disorder}).
  It is defined through a recurrence relation for
  the law of a random variable $\{R_n\}_{n=1,2, \ldots}$, which in absence of
  disorder ({\sl i.e.}, when the initial condition is degenerate) reduces to a
  particular case of the well-known Logistic Map. 
    The large-$n$ limit of the sequence of random variables $2^{-n}\log
  R_n$, a non-random quantity which is naturally interpreted as a free
  energy, plays a central role in our analysis. 
  The model 
  depends on a parameter $\alpha\in(0,1)$, related to the geometry of the
  hierarchical lattice, and has a phase transition in the sense that the free
  energy is positive if the expectation of $R_0$ is larger than a certain
  threshold value, and it is zero otherwise. It was conjectured in
  \cite{cf:DHV} that disorder is relevant (respectively, irrelevant or
  marginally relevant) if $1/2<\alpha<1$ (respectively, $\alpha<1/2$
  or $\alpha=1/2$), in the sense that an arbitrarily small amount of
  randomness in the initial condition modifies the critical point with
  respect to that of the pure ({\sl i.e.}, non-disordered) model if
  $\ga\ge1/2$, but not if $\ga<1/2$. Our main result  
  is a proof of  these conjectures for the case $\ga\ne 1/2$. We
  emphasize that for $\ga>1/2$ we find the {\sl correct scaling form}
  (for weak disorder) of the critical point shift. 
\end{abstract}

\maketitle

\section{Introduction}
\label{sec:intro}

\subsection{The model}
Consider the dynamical system 
defined by 
the initial condition $R_0^{(i)}>0$, $i\in \N :=\{1,2, \ldots\}$
and  the array of recurrence equations
\begin{eqnarray}
\label{eq:R}
  R_{n+1}^{(i)}\, =\, \frac{R_n^{(2i-1)}R_n^{(2i)}+(B-1)}B,  \ \ \  \, i \in \N,
\end{eqnarray}
for $n=0,1,\ldots$ and a given $B>2$. 
Of course if $R_0^{(i)}=r_0$ for every $i$, then the problem
 reduces to studying the quadratic recurrence equation
\begin{equation}
\label{eq:r}
r_{n+1}= \frac{r_n^2 + (B-1)}B,
\end{equation}
a particular case of a very classical problem, the {\sl logistic map},
as it is clear from the fact 
that $z_n:=1/2- r_n/(2(B-1))$ satisfies
the recursion 
\begin{equation}
\label{eq:logistic}
z_{n+1}=\frac{2(B-1)}B\,z_n(1-z_n).
\end{equation}
We are instead interested in non-constant initial data and, more
precisely, in initial data that are typical realizations of a sequence
of independent identically distributed (IID) random variables. In its
random version, the model was first considered in \cite{cf:DHV} (see
\S~\ref{sec:qmpm} and \S~\ref{sec:pinning} below for motivations in terms of
pinning/wetting models and for an informal discussion of what 
the
interesting questions are 
and what is expected to be true). We will
consider rather general distributions, but we will assume that all the
moments of $R^{(i)}_0$ are finite. As it will be clear later,
for our purposes it is actually useful to write
\begin{equation}
\label{eq:R0}
  R_0^{(i)} \, =\, \exp(\gb \go_i -\log \M (\gb)+h),
\end{equation}
with $\gb \ge 0$, $h \in \R$, 
 $\{\go_i\}_{i\in\N}$ a sequence of exponentially integrable IID
centered random variables normalized to  $\bbE\,\go_1^2=1$ and 
for every $\gb$
\begin{equation}
  \label{eq:M}
\M (\gb)\, := \, \bbE \exp(\gb \go_1)\, < \, \infty.
\end{equation}
The law of $\{\go_i\}_{i\in\N}$ is denoted by $\bbP$ and 
we will often alternatively denote the average $\bbE(\cdot)$
by brackets $\langle\cdot\rangle$.

Note that, for every $n$, $\{R_n^{(i)}\}_{i\in\N}$ are IID random
variables and therefore this dynamical system is naturally
re-interpreted as the evolution of the probability law $\cL _n$ (the
law of $R_n^{(1)}$): given $\cL_n$, the law $\cL_{n+1}$ is obtained by
constructing two IID variables distributed according to $\cL_n$ and
applying
\begin{equation}
\label{eq:basic}
R_{n+1}\, =\, \frac{R_n^{(1)}R_n^{(2)}+(B-1)}{B}.
\end{equation}

Of course, the iteration \eqref{eq:basic} is well defined for every
$B\ne0$. In particular, as detailed in Appendix
\ref{sec:Bsmallerthan2}, the case $B\in(1,2)$ can be mapped exactly
into the case $B>2$ we explicitly consider here, while for $B<1$ one
loses the direct statistical mechanics interpretation of the model
discussed in Section \ref{sec:pinning}.

\subsection{Quadratic maps and pinning models}
\label{sec:qmpm}

\begin{figure}[hlt]
\begin{center}
\leavevmode
\epsfxsize =14 cm
\psfragscanon
\psfrag{l0}{\small level $0$}
\psfrag{l1}{\small level $1$}
\psfrag{l2}{\small level $2$}
\psfrag{d0}{$d_0$}
\psfrag{d1}{$d_1$}
\psfrag{d2}{$d_2$}
\psfrag{d3}{$d_3$}
\psfrag{d4}{$d_4$}
\psfrag{u1}{$u_1$}
\psfrag{u2}{$u_2$}
\psfrag{u3}{$u_3$}
\psfrag{u4}{$u_4$}
\psfrag{traject}{\small trajectory $a$}
\psfrag{traject2}{\small trajectory $b$}
\epsfbox{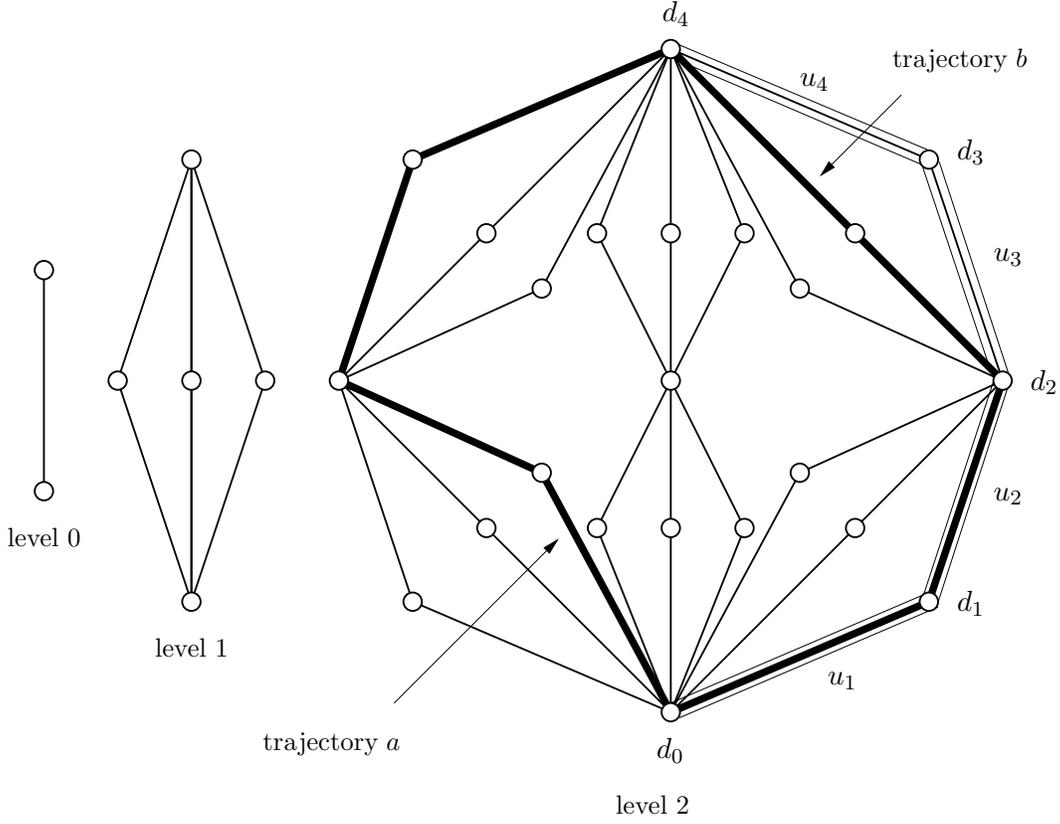}
\end{center}
\caption{Given $B=2,3, \ldots$ ($B=3$ in the drawing)
we build a diamond lattice by iterative steps
(left to right):
  at each step one replaces every bond by $B$
branches consisting of two bonds each. A trajectory of our process 
in a diamond lattice at {\sl level} $n$ is a path connecting the two {\sl poles}
$d_0$ and $d_{2^n}$: two trajectories, $a$ and $b$, are singled out by thick lines. 
Note that at level $n$, each trajectory is made of
$2^n$ bonds and there are $N_n$ trajectories, $N_0:=1$
and $N_{n+1}= B N_n^2$. A simple random walk at level $n$ is  the uniform measure
over the $N_n$ trajectories. 
A special trajectory, with vertices labeled 
$d_0, d_1, \ldots, d_{2^n}$, is  chosen (and marked by a triple line:
the right-most trajectory in the
drawing, but any other trajectory would lead to an equivalent model), we may call it defect line or wall boundary,
and rewards $u_j:= \gb \go_j - \log \M(\gb) +h$ (negative or positive) are assigned to the bonds of this trajectory.
The energy of a trajectory depends on how many and which bonds it shares
with the defect line: trajectory $a$ carries no energy, while trajectory $b$ carries energy $u_1+u_2$. 
The pinning model is then built by rewarding or penalizing 
the trajectories according to their energy  in the standard statistical mechanics fashion
and the partition function of such a model is therefore given by
$R_n$ in  \eqref{eq:DHV}. It is rather elementary, and fully detailed in   
 \cite{cf:DHV}, how
 to extract from \eqref{eq:DHV} the recursion
\eqref{eq:basic}. But the recursion itself is  well defined for
arbitrary real value 
$B\ne0$ and one may forget the definition of the
hierarchical lattice, as we do here. The definition of 
$\bP^B_n$ can also be easily generalized to $B >1$, see  \eqref{eq:EBgen} of
Appendix~\ref{sec:stand}.}
\label{fig:diamond}
\end{figure}

The model we are considering may be viewed as a hierarchical version
of a class of statistical mechanics models that go under the name of
(disordered) {\sl pinning} or {\sl wetting } models \cite{cf:Fisher,cf:Book}, that 
are going to be described in some detail in  \S~\ref{sec:pinning}.  It has been
introduced in \cite[Section 4.2]{cf:DHV}, where the partition function
$R_n=R_n^{(1)}$ is defined for $B=2,3, \ldots$ as
\begin{equation}
\label{eq:DHV}
R_n\, =\, \bE_n^B\left[\exp\left(\sum_{i=1}^{2^n} 
\left(\gb\go_i-\log \M (\gb)+h\right) \ind_{\left\{(S_{i-1},S_i)=(d_{i-1},d_i)\right\}}
  \right)\right],
\end{equation}
with $\{S_i\}_{i=0, \ldots , 2^n}$ a simple random walk (of law
$\bP_n^B$) on a hierarchical {\sl diamond}  lattice with growth parameter $B$
and $d_0, \ldots, d_{2^n}$ are the labels for the vertices of a particular path
that has been singled out and dubbed {\sl defect line}.
The construction of diamond lattices and a graphical description of the model are
 detailed in Figure~\ref{fig:diamond}
and its caption.

The phenomenon that one is trying to capture is the {\sl
  localization at (or delocalization away from) the defect line}, that is one
would like to understand whether the rewards (that could be negative,
hence penalizations) {\sl force} the trajectories to stick close to
the defect line, or the trajectories {\sl avoid} the defect line. A
priori it is not clear that there is necessarily a sharp distinction
between these two qualitative behaviors, but it turns out that it is
the case and which of the two scenarios prevails may be read from the
asymptotic behavior of $R_n$.  The Laplace asymptotics carries already
a substantial amount of information, so we define the {\sl quenched
  free energy}
\begin{equation}
  \label{eq:F}
\tf (\beta,h) \, := \, \lim_{n \to \infty}\frac1{2^n} \log R^{(1)}_n,
\end{equation}
where the limit is in the almost sure sense: the existence of such a
limit and the fact that it is non-random may be found in Theorem~\ref{th:F}.
Note in fact that $\partial_h \tf(\gb, h)$ coincides with
the $n \to \infty$ limit of 
$\bE^B_{n, \go} [ 2^{-n } \sum_{i}  \ind_{\left\{(S_{i-1},S_i)=(d_{i-1},d_i)\right\}}]$,
where $\bP^B_{n, \go}$ is the probability measure associated to the
partition function $R_n$, when $\partial_h \tf(\gb, h)$ exists
(that is for all $h$ except at most a countable number of points,
by convexity of $\tf (\gb, \cdot)$, see below). 
Therefore $\partial_h \tf(\gb, h)$ measures the density of contacts
between the walk and the defect line and below we will
see that $\partial_h \tf(\gb, h)$ is zero up to a critical value $h_c(\gb)$,
and positive for $h>h_c(\gb)$: this is a clear signature of a {\sl localization transition}.

\subsection{A first look at the role of disorder}
Of course if $\gb=0$ the {\sl disorder} $\go$ plays no role
and the model reduces to
  the one-dimensional map \eqref{eq:r}
 (in our language $\gb>0$
corresponds to the model in which  disorder is present).  This map
has two fixed points: $1$, which is stable, and $B-1$, which is
unstable.  More precisely, if $r_0<B-1$ then $r_n$ converges
monotonically (and exponentially fast) to $1$. If $r_0>B-1$, $r_n$
increases to infinity in a super-exponential fashion, namely
$2^{-n}\log r_n$ converges to a positive number which is of course
function of $r_0$.  The question is whether, and how,  introducing
disorder in the initial condition ($\gb>0$) modifies this behavior.

There is also an alternative way to link \eqref{eq:R} and \eqref{eq:r}. 
In fact, by  taking the average we obtain
\begin{equation}
\label{eq:Rav}
\Rav{n+1}\, =\, \frac{\Rav{n}^2 + (B-1)}{B},
\end{equation}
where we have dropped the superscript in $\langle R_n^{(i)}\rangle$.
Therefore the behavior of the sequence $\{ \Rav{n}\}_n$ is (rather)
explicit, in particular such a sequence tends (monotonically) to $1$
if $\Rav{0}<B-1$, while $\Rav{n}=B-1$ for $n=1,2, \ldots$ if
$\Rav{0}=B-1$.  This is already a strong piece of information on
$R_n^{(1)}$ (the sequence $\{ \cL_n\}_n$ is tight). Less informative
is instead the fact that $\Rav{n}$ diverges if $\Rav{0}>B-1$, even if
we know precisely the speed of divergence: in fact the sequence of random
variables can still be tight!  In principle such an issue may be tackled by looking at higher moments, but while
$\med{R_n}$ satisfies a closed recursion, the same is not true for
higher moments, in the sense that the recursions they satisfy depend
on the behavior of the lower-order moments.  For instance, if we set
$\gD_n := \text{var}\left( R_n \right)$, we have
\begin{equation}
\label{eq:Delta}
\gD_{n+1}\, =\, \frac{\gD_n \left( 2 \Rav{n}^2 + \gD_n\right)}{B^2}.
\end{equation}
In principle such an approach can be pushed further, but most important for 
understanding the behavior of the system is capturing the asymptotic behavior
of $\log R_n^{(i)}$, {\sl i.e.} \eqref{eq:F}.

\subsection{Quenched and annealed free energies}
Our first result says, in particular, that the quenched free energy
\eqref{eq:F} is well defined:

\medskip

\begin{theorem}
\label{th:F}
The limit in \eqref{eq:F} exists $\bbP(\dd \go)$-almost surely and in
$\mathbb L^1(\dd \mathbb P)$, it is almost-surely constant and it is
non-negative. 
The function $(\gb,h) \mapsto \tf(\gb , h+\log \M (\gb))$ is convex and 
$\tf(\beta,\cdot)$
   is non-decreasing (and convex). These properties are inherited from 
 $\tf_N(\cdot,\cdot)$, defined by
\begin{equation}
\tf_N(\gb,h)\,=\, \frac{1}{2^N}\med{\log R_N}.
\end{equation}
Moreover $\tf_N (\gb, h)$ converges to $\tf(\gb, h)$ with exponential speed, more precisely for
all $N\ge 1$
\begin{align}
  \tf_N(\gb,h)-2^{-N}\log B \le \tf(\gb,h)
  \le \tf_N(\gb,h)+2^{-N}\log \left(\frac{B^2+B-1}{B(B-1)}\right).
\label{eq:encadre}
\end{align}
\end{theorem}

\medskip

Let us also point out that $\tf(\beta,h)\ge 0 $ is immediate in view of the
fact that $R_n^{(i)}\ge (B-1)/B$ for $n\ge 1$, {\sl cf.} \eqref{eq:R}.
The lower bound $\tf(\beta,h)\ge 0 $ implies that we can split the
parameter space (or {\sl phase diagram}) of the system according to
$\tf(\gb,h)=0$ and $\tf (\gb, h)>0$ and this clearly corresponds to
sharply different asymptotic behaviors of $R_n$.  In conformity with
related literature, see \S~\ref{sec:pinning}, we define localized
and delocalized phases as $\cL:=\{(\gb, h):\, \tf(\gb,h)>0\}$ and
$\cD:=\{ (\gb, h):\, \tf(\gb,h)=0\}$ respectively.  It is therefore
natural to define, for given $\beta\ge0$, the {\sl critical} value
$h_c(\beta)$ as
\begin{equation}
h_c(\beta)\, =\, \sup\{h\in\R:\tf(\beta,h)=0\}.
\end{equation}
Theorem \ref{th:F} says in particular that 
\begin{equation}
h_c(\gb)\, =\, \inf\{h\in\R:\tf(\gb,h)>0\},
\end{equation}
and that $\tf(\gb,\cdot)$ is (strictly) increasing on $(h_c(\beta),\infty)$.
Note that, thanks to the properties we just mentioned, the contact fraction, defined in the end of \S~\ref{sec:qmpm},  is zero 
$h<h_c(\gb)$ and  is instead positive if
$h>h_c(\gb)$ (define the contact fraction by taking the inferior limit
for the values of $h$ at which  $\tf (\gb, \cdot)$ is not differentiable).

\smallskip

Another important observation on Theorem \ref{th:F} is that 
it yields also the existence of $\lim_{n \to \infty} 2^{-n} \log \langle R_n\rangle$
and this limit is simply $\tf(0,h)$, in  fact
$\tf_n(0, h)= 2^{-n }\log \Rav{n}$ for every $n$.
In statistical mechanics language $\Rav{n}$ is an {\sl annealed} quantity and
$\lim_{n \to \infty} 2^{-n} \log \langle R_n\rangle$
is the {\sl annealed free energy}: by Jensen
inequality it follows that $\tf(\beta,h)\le\tf(0,h)$
and $h_c(\beta)\ge h_c(0)$.
It is also a consequence of Jensen inequality (see Remark~\ref{rem:5}) 
the fact that $\tf(\gb, h+ \log \M( \gb))\ge \tf (0,  h)$,
so that $h_c(\beta)\le h_c(0)+ \log \M (\gb) $.
Summing up:
\begin{equation}
\label{eq:bb}
h_c(0) \, \le \, 
h_c(\beta)\, \le\, h_c(0)+ \log \M (\gb) .
\end{equation}
Therefore, by the convexity properties of $\tf(\cdot, \cdot)$ (Theorem~\ref{th:F})
and by \eqref{eq:bb}, we see that $h_c(\cdot) - \log \M(\cdot)$ is concave
and may diverge only at infinity, so that $h_c(\cdot)$ is a continuous function.

The following result on the {\sl annealed system}, {\sl i.e.} just the non-disordered system,
is going to play an important role:
 
\medskip

\begin{theorem}
  \label{th:pure} {\it (Annealed system estimates)}.  The function
  $h\mapsto \tf(0, h)$ is real analytic 
  except at $h=h_c:=h_c(0)$.  Moreover $h_c=\log
  (B-1)$ and there exists $c=c(B)>0$ such that for all
  $h\in(h_c,h_c+1)$
\begin{eqnarray}
\label{eq:alpha0}
  c(B)^{-1}(h-h_c)^{1/\alpha}\, \le\, \tf(0,h)\, \le
  \,  c(B)(h-h_c)^{1/\alpha},
\end{eqnarray}
where 
\begin{equation}
\label{eq:alpha}
  \alpha\, :=\, \frac{\log(2(B-1)/B)}{\log 2}.
\end{equation}
\end{theorem}

\medskip

Bounds on the annealed free energy can be extracted directly from
\eqref{eq:encadre}, namely
that for every $n\ge 1$
\begin{equation}
\label{eq:annbounds}
\frac{B(B-1)}{B^2+B-1}\exp\left(2^n \tf(0,h)\right) \, \le \, \Rav{n} \,\le\,
B \exp\left(2^n \tf(0,h)\right).
\end{equation}
Moreover let us note from now that $\alpha\in(0,1)$ and that
$1/\alpha>2$ if and only if $B<B_c:=2+\sqrt 2$, and $1/\alpha=2$ for
$B=B_c$.  It follows that $\tf(0, h)=o((h-h_c)^2)$ for $B<B_c$ ($\ga <1/2$),
while  this is not true for $B>B_c$ ($\ga>1/2$). 
\medskip

\begin{rem}\rm
\label{rem:oscillations}
  For  models defined on hierarchical lattices, in
  general
  one does not expect 
   the (singular part of the) free energy to have a pure
  power-law behavior close to the critical point $h_c$, but rather to
  behave like $H(\log (h-h_c))(h-h_c)^\nu$, with $\nu$ the critical
  exponent and $H(\cdot)$ a periodic function, see in particular \cite{cf:DIL}.  
    Note that, unless
  $H(\cdot)$ is trivial ({\sl i.e.} constant), the oscillations it produces become more and
  more rapid for $h\searrow h_c$. We have observed numerically such oscillations
  in our case and therefore we expect that
  estimate \eqref{eq:alpha0} cannot be improved at a qualitative level as $h$ approaches
  $h_c$
  (the problem of estimating sharply the size of the oscillations appears to be a
   non-trivial one, but this is not particularly important for our analysis).
\end{rem}

\subsection{Results for the disordered system}
\label{sec:results}

\smallskip

The first result we present gives  information
on the phase diagram: we use the definition
\begin{equation}
  \label{eq:D}
  \gD\,  = \, \gD(\gb)\, :=\, (B-1)^2\left(\frac{\M(2\gb)}{\M(\gb)^2}-1\right)\, \left(\ge\, 0\right),
\end{equation}
so that $  \text{Var}(R_0)\stackrel{h=h_c}= \gD$. The quantity $\gD$
should be though of  as the size of the disorder at a given $\gb$. 

\medskip

\begin{theorem}
\label{th:eps_c}
Recall that the critical value for the annealed system is $h_c= \log(B-1)$.
We have the following estimates on the quenched critical line:
\begin{enumerate}
\item
Choose $B \in (2, B_c)$. If  $\gD (\gb) \le  B^2 -2(B-1)^2$ then
$h_c(\gb)=h_c$. 
\item Choose $B>B_c$. Then $h_c(\gb)>h_c$ for every $\gb>0$.
  Moreover for $\gb$ small (say, $\gb \le 1$) one can find
  $c\in(0,1)$ such that
\begin{equation}
\label{eq:eps_c}
c \, \gb^{2\ga/(2\ga -1)} \, \le \, 
h_c(\gb ) -h_c \, \le \, c^{-1}  \gb^{2\ga/(2\ga -1)}.
\end{equation} 
\item If $B=B_c$ then one can find $C>0$ 
such that, for $\gb \le 1$,
\begin{equation}
\label{eq:shiftBc}
0 \le h_c(\gb)-h_c \le \exp(-C/\gb^2).  
\end{equation}
\end{enumerate}
Moreover if $\go_1$ is such that $\bbP(\go_1 >t)>0$ for every $t>0$,
then for every $B>2$ we have $h_c(\gb)-h_c>0$ for $\gb$ sufficiently
large, in fact $\lim_{\gb \to \infty} h_c(\gb )=\infty$.
\end{theorem}

Of course  \eqref{eq:shiftBc} leaves open an evident question  for
$B=B_c$, that will be discussed  in \S~\ref{sec:pinning}.  We point out
that the constant $C$ is explicit (see Proposition~\ref{th:lwbdc}) but
it does not have any particular meaning.  It is possible to show that
$C$ can be chosen arbitrarily close to the constant given in
\cite{cf:DHV}, but here, for the sake of simplicity, we have decided
to prove a weaker result ({\sl i.e.}, with a smaller constant).  This
is not a crucial issue, since the upper bound on $h_c(\gb)$ is not
comforted by a suitable lower bound.

\smallskip

\begin{figure}[h]
\begin{center}
\leavevmode
\epsfysize =5.7 cm
\psfragscanon
\psfrag{0}[c]{$0$}
\psfrag{D}[c]{$\gb$}
\psfrag{eps}[c]{$h$}
\psfrag{aa}[c]{ $\hat \gb$}
\psfrag{a}[c]{$\gb_c$}
\psfrag{Del}[c]{\large \cD}
\psfrag{Loc}[c]{\large \cL}
\psfrag{hc}[c]{$h_c$}
\psfrag{B>Bc}[c]{\large $B>B_c$}
\psfrag{B<Bc}[c]{\large $B<B_c$}
\epsfbox{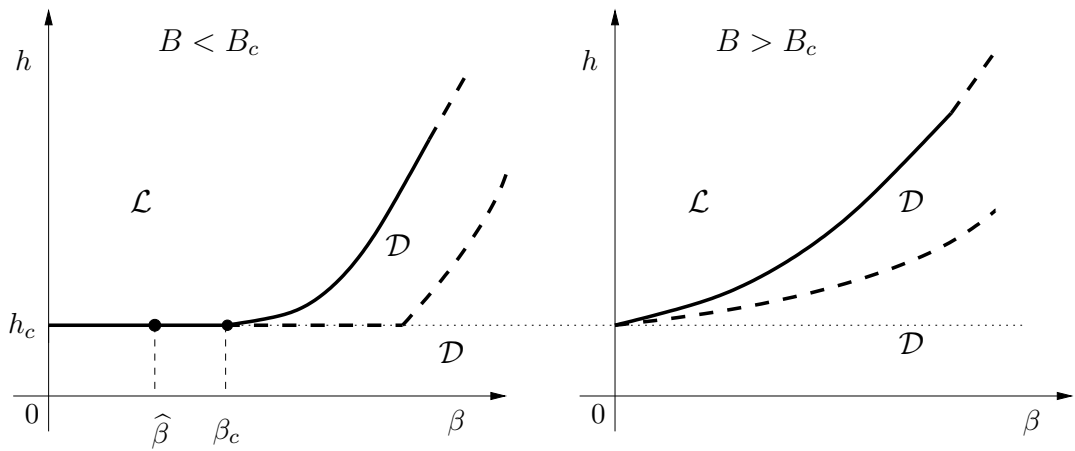}
\end{center}
\caption{\label{fig:hbeta} This is a sketch of the phase diagram 
 and a graphical view of
  Theorem~\ref{th:eps_c} and Theorem~\ref{th:paths}.  The thick line
  in both graphs is $h_c(\cdot)$. The dashed line is instead the lower bound
on $h_c(\cdot)$ which we obtain with our methods.
Below the dashed line we can establish the
  a.s. convergence of $R_n^{(i)}$ to $1$.  We have also used $\gb_c
  :=\sup\{\gb:\, h_c(\gb)=h_c\}$ and 
  $\widehat \gb := \sup\{\gb:\, \gD (\gb) < B^2-2(B-1)^2\}$. We do not prove  
  the (strict) inequality 
$\gb_c> \widehat \gb$.}
\end{figure}

The next result is about  the free energy.

\medskip

\begin{theorem}
\label{th:fe}
We have the following:
\begin{enumerate}
\item Choose $B\in (2,B_c)$ and $\gb$ such that
$\gD (\gb) < B^2-2(B-1)^2$. Then for every
$\eta\in (0,1)$ one can find $\epsilon>0$ such that
\begin{equation}
\label{eq:lbfe}
 \tf ( \gb, h)\, \ge\, (1-\eta)  \tf(0,h),
\end{equation}
for $h\in (h_c, h_c+\epsilon)$.
\item Choose $B>B_c$. Then for every
$\eta\in (0,1)$ one can find $c>0$ and $\gb_0>0$ such that
\eqref{eq:lbfe} holds for $\gb< \gb_0$ and 
$h-h_c \in (c \gb^{2\ga /(2\ga -1)}, 1)$.  
\end{enumerate}
\end{theorem}
\medskip

While the relevance of the analysis of the free energy will
be discussed in depth in the next subsection, it is natural to address 
the following issue: in a {\sl sharp} sense, how does the random 
array $R_n^{(1)}$ behave as $n$ tends to infinity?
We recall that the non-disordered system displays only three possible
asymptotic behaviors: $r_n\to 1$, $r_n=B-1$ for all $n$ and 
$r_n\nearrow \infty$ in a super-exponentially fast fashion.

What can be extracted directly from the free energy is quite
satisfactory if the free energy is positive: $R^{(1)}_n$ diverges at a
super-exponential speed that is determined to leading order.  However,
the information readily available from the fact that the free energy
is zero is rather poor; this can be considerably improved, starting
with the fact that, by the lower bound in \eqref{eq:encadre}, if the
free energy is zero then $\sup_n \langle \log R_n \rangle \le \log B$,
which implies the tightness of the sequence.

\medskip

\begin{theorem}
  \label{th:paths} If $\tf (\gb, h)=0$ then the sequence $\{R_n\}_n$
  is tight. Moreover 
if $h<h_c(\beta)$ then
\begin{equation}
\label{eq:paths}
\lim_{n \to \infty} R_n^{(1)}\, =\,1 \mbox{\; in probability}.
\end{equation}
\end{theorem}

\medskip

Let us mention that we also establish almost sure convergence of $R_n$
toward $1$ when we are able to find $\gamma\in (0,1)$ and $n\in \N$
such that $\bbE \left[ ([R_n-1]^+)^\gamma\right]$ is smaller than an
explicit constant (see Section~\ref{sec:rel}, in particular
Remark~\ref{rem:as1}). It is interesting to compare such results
with the estimates on the size of the partition
function $Z_{N,\go}$ of non-hierarchical pinning/wetting models, which
are proven in \cite[end of Sec.~3.1]{cf:T_fractmom} in the delocalized
phase, again via estimation of fractional moments of $Z_{N,\go}$
(which plays the role of our $R_n$).
\medskip 

What one should expect at criticality is rather unclear to us (see
however \cite{cf:GM_h} for a number of predictions and numerical
results on hierarchical pinning and
also \cite{cf:CEGM-JSP,cf:CEGM-CMP} for
some theoretical considerations on a different class of hierarchical models).

\subsection{Pinning models: the role of
disorder}
\label{sec:pinning}

Hierarchical models on diamond lattices, homogeneous or disordered
\cite{cf:Bleher,cf:BO,cf:CEGM-JSP,cf:CEGM-CMP,cf:DG}, are a powerful
tool in the study of the critical behavior of statistical mechanics
models, especially because real-space renormalization group
transformations {\sl \`a la } Migdal-Kadanoff are exact in this case.
In most of the cases, hierarchical models are introduced in association
with a more realistic non-hierarchical one. It should however be
pointed out that hierarchical models on diamond lattices are not rough
simplifications of non-hierarchical ones. They are in fact meant to
retain the essential features of the associated non-hierarchical
models (notably: the critical properties!). 
In particular, it would be definitely misleading to think of the
hierarchical model as a mean field approximation of the real one.

Non-hierarchical pinning models have an extended
literature ({\sl e.g.}  \cite{cf:Fisher,cf:Book}).  They may be
defined like in \eqref{eq:DHV}, with $S$ a symmetric random walk with
increment steps in $\{-1, 0, +1\}$, energetically rewarded or
penalized when the bond $(S_{n-1},S_n)$ lies on the horizontal axis
(that is $d_j=0$ for every $j$ in  \eqref{eq:DHV}),
but they can be restated in much greater generality by considering
arbitrary homogeneous Markov chains that visit a given site (say, the
origin) with positive probability and that are then rewarded or
penalized when passing by this site.  In their non-disordered version
\cite{cf:Fisher}, this general class of models has the remarkable
property of being {\sl exactly solvable}, while displaying a phase
transition --  a localization-delocalization transition -- and
the order of such a transition depends on a parameter of the model
(the tail decay exponent of the distribution of the first return of
the Markov chain to the origin: we call $\ga$ such an exponent and it
is the analog of the quantity $\ga$ in our hierarchical context, {\sl
  cf.} \eqref{eq:alpha}; one should however note that for
non-hierarchical models values $\ga\ge1$ can also be considered, in
contrast with the model we are studying here). As a matter of fact,
transitions of all order, from first order to infinite order, can be
observed in such models.  They therefore constitute an ideal set-up in
which to address the natural question: how does the disorder affect
the transition?

Such an issue has often been considered in the physical literature and
a criterion, proposed by A.~B.~Harris in a somewhat different context,
adapted to pinning models \cite{cf:FLNO,cf:DHV},  yields that the
disorder is irrelevant if $\gb$ is small and $\ga <1/2$, meaning by
this that quenched and annealed critical  points coincide and
the critical behavior of the free energy is the same for annealed and
quenched system (note that the annealed system is a homogeneous
pinning system, and therefore exactly solvable).  The disorder instead
becomes relevant when $\ga>1/2$, with a shift in the critical point
(quenched is different from annealed) and different critical behaviors
(possibly expecting a smoother transition, but the Harris criterion
does not really address such an issue). In the marginal case,
$\alpha=1/2$, disorder could be {\sl marginally} relevant or {\sl
  marginally} irrelevant, but this is an open issue in the physical
literature, see \cite{cf:FLNO,cf:DHV} and \cite{cf:Book} for further
literature.

\smallskip

Much progress has been made very recently in the mathematical literature 
on non-hierarchical pinning models, in particular:
\begin{enumerate}
\item The irrelevant disorder regime is under control
  \cite{cf:Ken,cf:T_cmp} and even more detailed results on the
  closeness between quenched and annealed models can be established
  \cite{cf:GTmuF}.
\item 
Concerning 
 the  relevant disorder regime, in \cite{cf:GT_cmp} it has been shown that the quenched
  free energy is smoother than the annealed free energy if $\ga>1/2$.
The non-coincidence of quenched and annealed critical points for 
{\sl large} disorder (and for every $\ga$) has been  proven in 
\cite{cf:T_fractmom} via an estimation of non-integer moments of the
partition function. The idea of considering non-integer moments (this time,
of $R_n-1$)  plays an important role also in the present paper.
\item A number of results on the behavior of the paths of the model
  have been proven addressing the question of what can be said about
  the trajectories of the system once we know that the free energy is
  zero (or positive) \cite{cf:GTdeloc,cf:GT_alea}.  One can in fact
  prove that if $\tf(\gb, h)>0$ then the process sticks close to the
  origin (in a strong  sense) and it is therefore in a localized
  ($\cL$) regime.  When  $\tf(\gb, h)=0$, and leaving aside the
  critical case, one expects that the process {\sl essentially never}
  visits the origin, and we say that we are in a delocalized regime
  ($\cD$).  We refer to \cite{cf:Book} for further discussion and
  literature on this point.
\end{enumerate}
\medskip

In this work we rigorously establish the full Harris criterion picture
for the hierarchical version of the model. In particular we wish to
emphasize that we do show that there is a shift in the critical point
of the system {\sl for arbitrarily small disorder} if $\ga>1/2$ and we
locate such a point in a window that has a precise scaling behavior,
{\sl cf.} \eqref{eq:eps_c} (a behavior which coincides with that
predicted in \cite{cf:DHV}).

As a side remark, one can also generalize the smoothing inequality
proven in \cite{cf:GT_cmp} to the hierarchical context and show that
for every $B>2$ there exists $c(B)<\infty$ such that, if $\go_1\sim
\cN (0,1)$, for every $\beta>0$ and $\delta>0$ one has
\begin{equation}
\label{eq:cB}
  \tf(\beta,h_c(\beta)+\delta)\, \le\,  {\delta^2}c(B)/{\beta^2},
\end{equation}
which implies that annealed and quenched free energy critical 
behaviors are different for $\ga >1/2$, {\sl cf.}  \eqref{eq:alpha0} (as
in \cite{cf:GT_cmp}, such inequality can be generalized well beyond
Gaussian $\go_1$, but we are not able to establish it only assuming
the finiteness of the exponential moments of $\go_1$). 
The proof of \eqref{eq:cB} is detailed in \cite{cf:LT}.

\medskip

Various intriguing issues remain open:
\begin{enumerate}
\item Is there a shift in the critical point at small disorder if
  $B=B_c$ (that is $\ga=1/2$)? We stress that in \cite{cf:DHV}
is predicted that $h_c(\gb)-h_c(0)\simeq \exp(-\log 2/\gb^2)$ for 
$\beta$ small.
\item Can one go beyond \eqref{eq:cB}? That is, can one find sharp estimates
  on the critical behavior when the disorder is relevant?
\item With reference to 
the caption of
  Figure~\ref{fig:hbeta}, can one prove $\gb_c > \widehat \gb$ (for 
  $B<B_c$)? 
\item Does the law of $R_n$ converge to a non-trivial limit for $n\to\infty$,
when $h=h_c(\gb)$?
 \end{enumerate}
\medskip

Of course, all these issues are open also in the non-hierarchical context
and, even if not every question becomes easier for the hierarchical model,
it may be the right context in which to attack them first.

\subsection{Some recurrent notation and organization of the subsequent sections} 

\label{sec:srn}

Aside for standard notation like $\lceil x\rceil := \min\{n\in \Z: \,
n\ge x\}$  and $\lfloor x\rfloor:= \lceil x\rceil -1$, or
$[\cdot]^+:=\max(0,\cdot)$, we will repeatedly use $\gD_n$ for the
variance of $R_n^{(1)}$, see \eqref{eq:Delta}, and $Q_n := \gD_n /
\Rav{n}^2$ so that from \eqref{eq:Rav} and \eqref{eq:Delta}, one sees
that
\begin{equation}
\label{eq:Q}
Q_{n+1}\, =\,
2 \left( \frac{B-1}B\right)^2
\left(\frac{\Rav{n}^4}{\Rav{n+1}^2 (B-1)^2}\right) 
\left( Q_n + \frac 12 Q_n^2\right),
\end{equation}
and we observe that
\begin{equation}
  \label{eq:Q0}
  Q_0\, =\,\left(\frac{M(2\gb)}{M(\gb)^2}-1\right)\stackrel{\beta\searrow0}
\sim\beta^2.
\end{equation}
Note that  $2(B-1)^2/B^2$ is smaller than $1$ if and
only if $B<B_c$ and 
\begin{equation}
\label{eq:Qbound}
\left(\frac{\Rav{n}^4}{\Rav{n+1}^2 (B-1)^2}\right) \, \le \, \left( \frac B{B-1}\right)^2.
\end{equation}

We will also frequently use  $P_n:=\langle R_n\rangle-(B-1)$, which satisfies
\begin{equation}
\label{eq:P}
P_{n+1}\, =\,
2\frac{(B-1)}B P_n+ \frac 1B P_n^2,
\end{equation}
and $P_0=\gep$ in our notations (see \eqref{eq:gep} below).
With some effort, one can explicitly verify that for every $n$ 
\begin{equation}
\label{eq:QboundP}
\left(\frac{\Rav{n}^4}{\Rav{n+1}^2 (B-1)^2}\right) \, \le \, 1+ \frac{4P_n}{B(B-1)}.
\end{equation}

Finally, there is some notational convenience at times in making the change of variables 
\begin{equation}
  \label{eq:gep}
  \gep\, :=\, \Rav{0}-(B-1)\, =\,  e^h-(B-1),
\end{equation}
and 
\begin{equation}
\label{eq:hatF}
\hat
\tf(\gb,\gep)\,:=\,\tf(\beta,h(\gep)),
\end{equation}
and when we write $ h(\gep)$ we refer to the invertible map 
defined by \eqref{eq:gep}.

\medskip

The work is organized as follows. Part (1) of Theorem \ref{th:eps_c}
and of Theorem \ref{th:fe} are proven in Section \ref{sec:irrel}. In
Section \ref{sec:LBrel} we prove part (2) of Theorem \ref{th:fe} and,
as a consequence, part (2) of Theorem \ref{th:eps_c}, except the lower
bound in \eqref{eq:eps_c}. Part (3) of Theorem \ref{th:eps_c} is
proven in Section \ref{sec:LBmarg} and the lower bound of
\eqref{eq:eps_c} in Section \ref{sec:rel} (after a brief sketch of our
method). The proof of Theorem \ref{th:paths} is given in Section
\ref{sec:deloc}. Finally, the proofs of Theorems \ref{th:F} and
\ref{th:pure} are based on more standard techniques and can be found
in Appendix \ref{sec:stand}.

\section{Free energy lower bounds: $B<B_c= 2+\sqrt{2}$}

\label{sec:irrel}

We want to give a proof of part (1) of Theorem \ref{th:fe}, which in particular
implies part (1) of Theorem \ref{th:eps_c}. 

The strategy goes roughly as follows: since $h>h_c$ is close to $h_c$,
that is $\gep (=P_0)>0$ is close to $0$, $P_n$ keeps close to zero for
many values of $n$ and $P_{n+1} \approx (2(B-1)/B)P_n$ (recall \eqref{eq:P} and the fact that
$2(B-1)/B>1$ for $B>2$).  This is going to be true up to $n$ much
smaller than $\log (1/\gep) / \log (2(B-1)/B)$.  At the same time for
the normalized variance $Q_n$ we have the approximated recursion
$Q_{n+1} \approx 2((B-1)/B)^2 (Q_n+(1/2) Q_n^2)$, which one derives
from \eqref{eq:Q} by using $P_n\approx 0$. Since $2((B-1)/B)^2 <1$ is
equivalent to $B<B_c$, we easily see that (if $Q_0$ is not too large)
$Q_n$ shrinks at an exponential rate. This {\sl scenario} actually
breaks down when $P_n$ is no longer small, but at that stage $Q_n$ is
already extremely small (such a value of $n$ is precisely defined and
called $n_0$ below). From that point onward $Q_n$ starts growing
exponentially and eventually it diverges, but after $(1+\gamma) n_0$
steps, for some $\gamma>0$, $Q_n$ is still small while $P_n$ is large,
so that a second moment argument, combined with \eqref{eq:encadre}
which yields a control on $\tf( \gb, h)$ via $\tf_n(\gb, h)$, allows
to conclude.




\medskip

Before starting the proof we give an upper bound on the size 
of  $Q_n(= \gD_n/ \Rav{n}^2)$ in the regime in which the recursion for
$\Rav{n}$ can be linearized (for what follows, recall  
\eqref{eq:Q},
\eqref{eq:Q0} and \eqref{eq:Qbound}).

\smallskip

\begin{lemma}
\label{th:linear}
Let $B\in(2,B_c)$ and $\gb$ such that $\gD=\gD (\gb)<B^2-2(B-1)^2$.  There exist
$c:=c(B,\Delta)>0$, $c_1:=c_1(B,\Delta)>0$ and
$\gd_0:=\delta_0(B,\Delta)>0$ with
\begin{equation}
\label{eq:d0c1}
  2(1+\delta_0)\left(\frac{B-1}B\right)^2\, <\,1,
\end{equation}
such that for every $\gep$ satisfying $0<  \gep/(B-1) < 
( (B^2-2(B-1)^2)/\gD)^{1/2}-1
$ (recall the definition \eqref{eq:gep} of $\gep$)
and
\begin{equation}
\label{eq:x0-1}
n \, \le \,n_0:=  \left\lfloor \log\left(c\,\gd_0 /\gep \right) / 
\log \left( \frac{2(B-1)}B\right)\right\rfloor,
\end{equation}
one has
\begin{equation}
\label{eq:linear}
Q_n \, \le \, c_1\, \left( 2(1+\gd_0)\left(\frac{B-1}B\right)^2 \right)^n\,Q_0 .
\end{equation} 
\end{lemma}
\smallskip

Note that the condition 
on $\gep $ simply guarantees  $\gD_0 =
(1+\gep/(B-1))^2 \gD$
is smaller than
$ B^2-2(B-1)^2$.
\smallskip 
 
\noindent
{\it Proof of Lemma \ref{th:linear}.}
Recall that  $P_n=\langle R_n\rangle-(B-1)$ and that
it satisfies the recursion \eqref{eq:P}
(and that $P_0=\gep$).

 For $G_n := (P_n/P_0) (2(B-1)/B)^{-n}$ we have from \eqref{eq:P} and
\eqref{eq:gep}
\begin{equation}
G_{n+1} \, =\, G_n  + \frac {\gep}B \,\left( 2\frac{(B-1)}B\right)^{n-1} G_n^2,
\end{equation}
and $G_0=1$. If  $G_m\le 2$ for $m\le n$, then
\begin{equation}
\frac{G_{n+1}}{G_n} \, \le \, 1+ 2\frac {\gep}B \,\left( 2\frac{(B-1)}B\right)^{n-1},
\end{equation}
which entails

\begin{equation} 
\label{eq:recineqGn}
G_{n+1}  \, \le \, 
\exp\left( 2\frac {\gep}B \sum_{j=0}^{n}\left( 2\frac{(B-1)}B\right)^{j-1}
\right)
\, \le \, 
1+\gep C(B) \left(\frac{2(B-1)}B\right)^{n+1},
\end{equation}
for a suitable constant $C(B)<\infty$. 

As we have already remarked, our assumption on $\gep$ yields  $\gD_0 < B^2-2(B-1)^2$,
so  
\begin{equation}
  \label{eq:cQ0}
Q_0 \, <\,  \left(\frac B{B-1}\right)^2-2. 
\end{equation}
Choose $\delta_0>0$ sufficiently small so that \eqref{eq:d0c1} is
satisfied and moreover
\begin{equation}
  \label{eq:d0c2}
  2\left(\frac{B-1}B\right)^2(1+\delta_0)(Q_0+\frac12Q_0^2)\, <\,Q_0, 
\end{equation}
(the latter can be satisfied in view of \eqref{eq:cQ0}).  It is
immediate to deduce from \eqref{eq:recineqGn} that if $c$ in 
\eqref{eq:x0-1} is chosen sufficiently small (in particular, $c\le B(B-1)/8$), 
then $G_n\le 2$ for
$n\le n_0$ and, as an immediate consequence,
\begin{equation}
\label{eq:recineqPn}
0<P_{n} \, \le \, 2 \gep  \left(\frac{2(B-1)}B\right)^{n}\, \le \,
2c\delta_0\,\le\, 
\delta_0 \frac {B(B-1)}4,
\end{equation}
where the first inequality is immediate from \eqref{eq:P} and $P_0=\gep>0$.
Now we apply \eqref{eq:QboundP}
\begin{eqnarray}
  \label{eq:Qnnuova}
  Q_{n+1}\le 2\left(\frac{B-1}B\right)^2(1+\delta_0)(Q_n+\frac12Q_n^2).
\end{eqnarray}
Notice also that $Q_1<Q_0$ thanks to \eqref{eq:d0c2}. From this it is
easy to deduce that, as long as $n\le n_0$, $Q_n$ is decreasing and
satisfies \eqref{eq:linear} for a suitable $c_1$.  In particular, $c_1(B, \gD_0)$ can be chosen such that
$\lim_{\gD_0 \searrow 0} c_1(B, \gD_0)=1$.  \qed

\bigskip

\noindent
{\it Proof of Theorem \ref{th:fe}, part (1).}
We use the bound \eqref{eq:Qbound} to get
\begin{equation}
\label{eq:superlinear}
Q_{n+1} \, \le \, 2 \left( Q_n + \frac 12 Q_n^2 \right) \, \le \, 3Q_n,
\end{equation}
where the last inequality holds as long as $Q_n \le 1$.  
Then we apply Lemma~\ref{th:linear} (recall in particular $\gd_0$ and $n_0$ in there).
Combining
\eqref{eq:linear} and \eqref{eq:superlinear} we get
\begin{equation}
\label{eq:QnQn0}
Q_n \, \le \, Q_{n_0} 3^{n-n_0} \, \le\, 
c_1\,Q_0 \left( 2(1+\gd_0)\left(\frac{B-1}B\right)^2 \right)^{n_0} 3^{n-n_0},
\end{equation}
for every $n\ge n_0$ 
satisfying $Q_n \le 1$ (which implies $Q_n'\le 1$ for all $n'\le n$ as $Q_n$ is increasing). Of course this boils
down to requiring that the right-most term in \eqref{eq:QnQn0} does not
get larger than $1$. Since $n_0$ diverges as $\gep \searrow 0$, if we choose 
$\gamma >0$ such that $3^\gamma\, 2(1+\gd_0) (B-1)^2/B^2 <1$, then 
the right-most term in \eqref{eq:QnQn0} is bounded above for every
$n \le (1+ \gamma) n_0$ by a quantity $o_\gep(1)$ which 
vanishes for $\gep\to0$.
Summing all up:
\begin{equation}
Q_{\lfloor(1+\gamma)n_0\rfloor }\, =\, o_\gep(1).
\end{equation}
 Next,
note that 
\begin{equation}
\label{eq:decomp}
\langle \log R_{\lfloor(1+\gamma)n_0\rfloor}\rangle  \ge
   \log\left(\frac 12
  \med{R_{\lfloor(1+\gamma)n_0\rfloor}}\right)
\bbP\left(R_{\lfloor(1+\gamma)n_0\rfloor}\ge
  \frac 12 \med{R_{\lfloor(1+\gamma)n_0\rfloor}}\right)+\log
  \left(\frac{B-1}B\right),
\end{equation}
where we have used the fact that $R_n\ge
(B-1)/B$ for $n \ge 1$. Applying the Chebyshev inequality one has
\begin{eqnarray}
\label{eq:PZ}
 \bbP\left(R_{\lfloor(1+\gamma)n_0\rfloor}\ge (1/2)
\med{R_{\lfloor(1+\gamma)n_0\rfloor}}\right)&\ge&
1-4Q_{\lfloor(1+\gamma)n_0\rfloor}=1+o_\gep(1).
\end{eqnarray}
Therefore, from  \eqref{eq:annbounds}, \eqref{eq:decomp} and
\eqref{eq:PZ}  one has
\begin{equation}
  \tf_{\lfloor(1+\gamma)n_0\rfloor}(\beta,h)\,\ge\,
(1+o_\gep(1)) \,\hat \tf(0,\gep)-
2^{-\lfloor(1+\gamma)n_0\rfloor}c(B),
\end{equation}
for some $c(B)<\infty$ and, from \eqref{eq:encadre} (or, equivalently, \eqref{eq:nonincrease}),
\begin{eqnarray}
 \hat \tf(\beta,\gep )\ge (1+o_\gep(1)) \,\hat\tf(0,\gep)-
2^{-\lfloor(1+\gamma)n_0\rfloor}c_1(B).
\end{eqnarray}
Since $\hat\tf(0,\gep)2^{\lfloor(1+\gamma)n_0\rfloor}$ diverges for
$\gep\to0$ if $\gamma>0$, as one may immediately check from
\eqref{eq:x0-1} and \eqref{eq:alpha0}, one directly extracts that for
every $\eta>0$ there exists $\gep_0>0$ such that
\begin{equation}
\label{eq:R1alpha}
\tf(\beta,h) \, =\, 
\hat \tf(\beta,\gep)\,
\ge \, (1-\eta) \hat\tf(0,\gep)\, = \, (1-\eta) \tf(0,h),
\end{equation}
for $\gep \le \gep_0$, {\sl i.e.} $h \le h_c(0) + \log(1+\gep/(B-1))$, and we are done.
\qed 

\section{Free energy lower bounds: $B\ge B_c =2+\sqrt{2}$}

\label{sec:LBrel}

The arguments in this section are close in spirit to the ones of the
previous section. However, since $B>B_c$, the constant $2((B-1)/B)^2$ in
the linear term of the recursion equation \eqref{eq:Q} is larger than
one, so the normalized variance $Q_n$ grows from the very beginning.
Nonetheless, if $Q_0$ is small, it will keep small for a while. The
point is to show that, if $P_0$ is not too small (this concept is of
course related to the size of $Q_0$), when $Q_n$ becomes of order one
$P_n$ is sufficiently large. Therefore, once again, a second moment
argument and \eqref{eq:encadre} yield the result we are after, that
is:

\medskip

\begin{proposition}
\label{th:LB_B>B_c}
Let $B>B_c$.  For every $\eta \in (0,1)$ there exist $c>0$
and $ \gb_0>0$ such that
\begin{equation}
 \label{eq:LBBc}
\tf \left( \beta ,h\right) \, \ge \, 
(1-\eta)  \tf \left( 0, h\right),
\end{equation}
for $\gb \le \gb_0$ and $c \gb^{2\ga/ (2\ga -1)}\le h-h_c(0)\le 1$.
This implies in particular that
$h_c(\gb) < h_c(0)+ c \gb^{2\ga/ (2\ga -1)}$, for every $\gb \le \gb_0$.
\end{proposition}
Of course this proves part (2) of Theorem \ref{th:fe} and the upper bound
in \eqref{eq:eps_c}.

\medskip

In this section $q:=2(B-1)^2/B^2$ and $\bar q := 2(B-1)/B$: note that
in full generality $ q < \bar q <2$ and $\bar q>1$, while $q>1$
because we assume $B>B_c$. One can easily check that
\begin{eqnarray}
  \label{eq:aa-1}
\frac{\ga}{2\ga-1}=\frac{\log\bar q}{\log q}.
\end{eqnarray}
Moreover in what follows
some expressions are in the form $\max A$, $A\subset \N \cup\{0\}$:
also when we do not state it explicitly, we do assume that $A$ is not
empty (in all cases this boils down to choosing $\gb$ sufficiently
small).

\smallskip

We start with an upper bound on the growth of 
$\Rav{n}= (B-1)+ P_n$ (recall \eqref{eq:P}) 
for $n$ {\sl not too large}.
\medskip

\begin{lemma}
\label{th:lem1-rel}
If $P_0=c_1 \gb^{2\ga/(2\ga -1)}$, $c_1>0$, then 
\begin{equation}
\label{eq:lem1-rel}
P_n \,\le\, 2c_1 \gb^{2\ga/(2\ga -1)} \bar q ^n\, \le \, 1, 
\end{equation}
for $n \le N_1 := \max\{n:\, C_1(B) c_1 \gb^{2\ga/(2\ga -1)} \bar q ^n
\le 1 \}$, where 
\begin{equation}
C_1(B) \,:= \,
2 \max\left(\frac1{(\bar q -1)B \log 2},1\right).
\end{equation}
\end{lemma}

\medskip

The next result controls the growth of the variance of $R_n$ in the
regime when $\Rav{n}$ is close to $(B-1)$, {\sl i.e.} $P_n$ is small.
Let us set $N_2 := \max\{n:\, (2 c_1 /(\bar q -1))\gb ^{2\ga
  /(2\ga -1)} \bar q^n \le (\log 2)/2\}$. Observe that $N_2\le N_1$
and recall that $Q_0 \stackrel{\gb \searrow 0} \sim \gb^2$, cf.
\eqref{eq:Q0}.

\medskip
\begin{lemma}
\label{th:lem2-rel}
Under the same assumptions as in  Lemma~\ref{th:lem1-rel}, for
$Q_0 \le 2\gb^2$ and assuming $c_1 \ge 20^{\log \bar q / \log q}$
we have 
\begin{equation}
Q_n \, \le\, 2 Q_0 q^n,
\end{equation}
for $n \le N_2$.
\end{lemma}

\medskip

\noindent
{\it Proof of Proposition~\ref{th:LB_B>B_c}}.
Let us choose $c_1$ as in  Lemma~\ref{th:lem2-rel}. Let us observe also that,
thanks to \eqref{eq:aa-1},
 $N_2 = \lfloor \log (1/\gb^2) /\log q - \log (C c_1)/ \log \bar q\rfloor$ 
 for a suitable choice of the constant $C=C(B)$. 
Therefore   Lemma~\ref{th:lem2-rel} ensures that
\begin{equation}
Q_{N_2} \, \le\, 4 (C c_1)^{-\log q /\log \bar q}.
\end{equation}
From the definition of $Q_n$ we directly see that
$Q_{n+1 }\le 3 Q_n$ if $Q_n \le 1$, as in \eqref{eq:superlinear}. 
Therefore for any fixed $\gd\in(0,1/16)$
\begin{equation}
\label{eq:Qbound1}
Q_{N_2+n } \,\le \,3^n 4(C c_1)^{-\log q /\log \bar q} \, \le \, 4\gd,
\end{equation}
if 
\begin{equation}
\label{eq:Qbound2}
n \, \le \, N_3 \, :=\, \left \lfloor \frac{\log q \log (C c_1)}{\log
  \bar q \log 3} - \frac{\log(1/\gd)}{\log 3}\right\rfloor.
\end{equation}
Since $Q_{N_2+N_3} \le
4\gd$ (by definition of $N_3$), we have then
 \begin{equation}
   \bbP\left( R_{N_2+N_3} \, \le \, \frac 12 \Rav{N_2+N_3} \right) \, 
\le \, 16 \gd.
 \end{equation}
As a consequence, applying \eqref{eq:encadre} and \eqref{eq:annbounds}
with $N=N_2+N_3$ one finds
\begin{eqnarray}
  \label{eq:c3}
  \tf(\beta,h)\ge (1-16\delta)\tf(0,h)-2^{-(N_2+N_3)}c_3(B),
\end{eqnarray}
of course with $h$ such that $P_0=c_1\beta^{2\ga/(2\ga-1)}$, {\sl i.e.},
\begin{eqnarray}
  h=\log\left((B-1)+c_1\beta^{2\ga/(2\ga-1)}\right).
\end{eqnarray}
  The last
step consists in showing that the last term in the right-hand side of
\eqref{eq:c3} is negligible with respect to the first one.  A look at
\eqref{eq:Qbound2} shows that $N_3$ can be made arbitrarily large by
choosing $c_1$ large; moreover, by definition of $N_2$ we have
\begin{equation}
2^{N_2} c_1^{1/\ga}\gb ^{2/(2\ga -1)} \, \ge \, \frac 12 C^{-1/\ga},
\end{equation}
for $\gb $ sufficiently small.  From these two facts and from the
critical behavior of $\tf(0,\cdot)$ (cf. \eqref{eq:alpha0}) one
deduces that for any given $\delta$ one may take $c_1$ sufficiently
large so that
\begin{equation}
2^{-(N_2+N_3)}/\tf(0,h)\le \delta,
\end{equation}
provided that $h\le h_c(0)+1$.
For a given $\eta\in(0,1)$ this  proves \eqref{eq:LBBc} whenever 
$\beta$ is sufficiently small and $c \gb^{2\ga/(2\ga-1)}\le h-h_c(0)\le1$,
with $c$ sufficiently large (when $\eta$ is small) but independent of $\beta$.


\qed

\medskip

\noindent
{\it Proof of Lemma~\ref{th:lem1-rel}.}
Call $N_0$ the largest value of $n$ for which $P_n \le 2c_1 \gb^{2\ga/(2\ga -1)}
\bar q^n$ (for $c_1$ and $\gb$ such that $P_0 \le  1$).
Recalling \eqref{eq:P}, for $n \le N_0$ we have
\begin{equation}
\frac{P_{n+1}}{P_n}\, \le \, \bar q \left( 1+ \frac{2c_1}{B \bar q}
\gb ^{2\ga/(2\ga -1)} \bar q ^n\right),
\end{equation}
so that for $N\le N_0$, using the properties  of $\exp(\cdot)$ and
the elementary bound $\sum_{n=0}^{N-1} a^n\le a^N/(a-1)$ ($a>1$),
 we obtain
\begin{equation}
P_N \, \le \, P_0 \,\bar q^N\exp \left( \frac{2c_1}{ (\bar q-1)B}
\gb ^{2\ga/(2\ga -1)} \bar q ^N \right).
\end{equation}
The latter estimate yields a lower bound on $N_0$:
\begin{equation}
\label{eq:ifrhs}
N_0 \, \ge \, \, \max\left\{
n:\, \frac{2c_1}{ (\bar q-1)B}
\gb ^{2\ga/(2\ga -1)} \bar q ^n \le \log 2
\right\}.
\end{equation}
$N_1$ is found by choosing it as the minimum between the right-hand side 
in \eqref{eq:ifrhs} and the maximal value of $n$ for which 
the second inequality in \eqref{eq:lem1-rel} holds.
  \qed

\medskip

\noindent
{\it Proof of Lemma~\ref{th:lem2-rel}.}
Let us call $N_0^{\prime}$ the largest $n$ such that $Q_n \le 2Q_0 q^n$ 
($N_0^{\prime}$ is introduced to control the nonlinearity in \eqref{eq:Q}) 
and let us work with $n \le \min(N_0^{\prime}, N_2)$.
Since $N_2 \le N_1$, 
 ($N_1$ given in Lemma~\ref{th:lem1-rel}), the bound \eqref{eq:lem1-rel} holds and $P_n \le 1$.
Therefore, by using first \eqref{eq:Q} and \eqref{eq:QboundP}, and then   
 \eqref{eq:lem1-rel}, we have 
\begin{equation}
\label{eq:satQ}
\frac{Q_{n+1}}{Q_n} \, \le \, q (1+P_n) \left( 1+2\gb^2 q^n\right)\, \le \, 
q \left(1 + 2 c_1  \gb ^{2\ga /(2\ga -1)} \bar q^n + 
4\gb^2 q^n
\right), 
\end{equation}
which implies
\begin{equation}
\label{eq:satQ-2}
Q_n \, \le \, 
Q_0 q^n \exp \left( 
  \frac{2c_1}{\bar q-1} \gb ^{2\ga /(2\ga -1)}\bar q^n  + \frac{4 }{q-1} \gb^2 q^n
\right). 
\end{equation}
By definition of $N_2$ the first term in the exponent is at most $(\log 2)/2$.
Moreover $n \le N_2$ implies, via \eqref{eq:aa-1},
\begin{equation}
\label{eq:satQ-3}
n \, \le \, \frac{\log (1/\gb^2)}{\log q}- \frac{\log\left((4/\log 2) c_1
/ (\bar q -1)\right)}{ \log \bar q},
\end{equation}
and one directly sees that for such values of $n$ we have $\gb^2 q^n
\le \left(4 c_1 / (\bar q -1)\log 2\right)^{-\log q / \log \bar
q}$. Therefore also the second term in the exponent ({\sl cf.}
\eqref{eq:satQ-2}) can be made smaller than $(\log 2)/2$ by choosing
$c_1$ larger than a number that depends only on $B$, see the statement
for an explicit expression.

Summing all up, for $c_1$ chosen suitably large, $Q_n\le 2Q_0 q^n$ for
$n \le \min(N_0^{\prime}, N_2)$. But, by definition of $N_0^{\prime}$,
this just means $n \le N_2$ and the proof is complete.  \qed

\subsection{The $B=B_c$ case}

\label{sec:LBmarg}
\begin{proposition}\label{th:lwbdc}
Set $B=B_c$. There exists $\gb_0$ such that for all $\gb\le \gb_0$
\begin{equation}
h_c(\gb)-h_c(0)\, < \,  \exp\left(-\frac{(\log2)^2}{2 \gb^2}\right).
\end{equation}
\end{proposition}

\medskip

\begin{rem} \rm The constant $(\log 2)^2/2$ that appears in the exponential is
  certainly not the best possible. In fact, one can get arbitrarily close 
   to the optimal constant $\log 2$ given in \cite{cf:DHV},
  but we made the choice to keep the proof as simple as possible.
\end{rem} 

\medskip

\noindent{\it Proof of Proposition \ref{th:lwbdc}.}
Choose 
\begin{equation}
h\,= \, e^{-(\log 2)^2/ (2\gb^2)} + \log(B_c-1),
\end{equation}
 so that 
 \begin{equation}
 P_0\, =\, 
\exp(h)-(B_c-1) \stackrel{\gb \searrow 0}\sim (B_c-1) \exp(-(\log 2)^2/(2
\gb^2)).
\end{equation}
 Given $\gd>0$ small (for example, $\gd=1/70$), we let
$n_\gd$ be the integer uniquely identified (because of the strict
monotonicity of $\{P_n\}_n$) by
\begin{equation}
P_{n_\gd}\, < \, \gd \, \le \, P_{n_\gd +1},
\end{equation} 
(we assume that $P_0<\gd$, which just means that we take $\gb$ small
enough).  We observe that \eqref{eq:P} implies $P_{n+1}/P_n \ge
\sqrt{2}$ for every $n$, from which follows immediately that (say, for
$\gb$ sufficiently small)
\begin{eqnarray}
n_\gd\le  \left \lceil \frac{\log 2}{\beta^2}\right\rceil.
\end{eqnarray}
  We want to show first of
all that $Q_{n_\gd}$ is of the same order of magnitude as
$Q_0$, and therefore much smaller than $P_{n_\gd}$ (for $\gb$ small)
in view of $Q_0 \stackrel{\gb \searrow 0}\sim \gb^2$.

From \eqref{eq:Q}, recalling the definition of $P_n$ ({\sl cf.}
\eqref{eq:P}) and the bound \eqref{eq:QboundP}, we derive
\begin{equation}
 \label{eq:QQ}
Q_{n+1}\, =\,
\left(\frac{\Rav{n}^4}{\Rav{n+1}^2 (B-1)^2}\right) 
\left( Q_n + \frac 12 Q_n^2\right)  
\, \le\,Q_n\left (1+P_n\right)\left(1+\frac{Q_n}{2}\right).
\end{equation}
If we define $c(\gd)$ through
\begin{eqnarray}
  \label{eq:cd}
  c(\gd)=\prod_{k=0}^{\infty}\left(1+\gd 2^{-k/2}\right)\le 
\exp(\delta (2+\sqrt 2)) \, \le \, \frac{21}{20},
\end{eqnarray}
from \eqref{eq:QQ} we directly obtain that,
as long as $Q_n\le 3 Q_0$ and $n\le n_\gd$,
\begin{eqnarray}
  Q_n\le Q_0 \left(1+(3/2)Q_0\right)^n\,\prod_{k=0}^{n-1}(1+P_n)\le
 c(\gd)\,Q_0\, e^{(3/2)Q_0\,n}.
\end{eqnarray}
It is then immediate to check, using \eqref{eq:Q0}, that $Q_{n_\gd}\le 3 Q_0$ for $\gb $ small.

But, as already exploited in \eqref{eq:superlinear},
$Q_{n+1}/Q_n \le 3$ for every $n$ such that $Q_n \le 1$,
so that  $Q_{n_\gd +n } \le  4\gb^2 3^{n} \le 1$ for 
$n\le n_1:=  \log_3 (1/(4\gb^2))-1 $.
But for such values of $n$
\begin{equation}
P_{n_\gd+n} \, \ge \,  \gd 2^{(n-1)/2},
\end{equation}
so that we directly see that $P_{n_\gd+n_1}$ diverges as $\gb$
tends to zero, and therefore $\Rav{n_\gd+n_1}$, can be made large for
$\gb$ small, while $Q_{n_\gd+n_1}$, that is the ratio between the
variance of $R_{n_\gd+n_1}$ and $\Rav{n_\gd+n_1}^2$ is bounded by $1$.  By
exploiting $R_n \ge (B-1)/B$ for $n \ge 1$ and using Chebyshev
inequality it is now straightforward to see that $\langle \log
(R_{n_\gd+n_1}/B)\rangle >0$ and by \eqref{eq:encadre} (or,
equivalently, \eqref{eq:nonincrease}) we have $\tf (\gb, h)>0$.

\qed

\section{Free energy upper bounds beyond annealing}

\label{sec:rel}

In this section we introduce our main new idea, which we briefly
sketch here.  In order to show that the free energy vanishes for $h$
larger than, but 
close to, $h_c(0)$, we take the system at the
$n$-th step of the iteration, for some $n=n(\gb)$ that scales suitably
with $\gb$ (in particular, $n(\gb)$ diverges for $\gb\to0$) and we
modify (via a tilting) the distribution $\bbP$ of the disorder.  If
$\alpha>1/2$, it turns out that one can perform such
 tilting so to guarantee on one hand that, under the new law, $R_{n(\gb)}$
is concentrated around $1$, and, on the other hand,  that the two laws are
very close (they have a mutual density close to $1$).  This in turn
implies that $R_{n(\gb)}$ is concentrated around $1$ also under the
original law $\bbP$, and the conclusion that $\tf(\gb,h)=0$ follows
then via the fact that if  some non-integer moment (of order
smaller than $1$) of $R_{n_0}-1$ is sufficiently small for some
integer $n_0$, then it remains so for every $n\ge n_0$ (cf.
Proposition \ref{th:fractmom}).

\subsection{Fractional moment bounds}
The following result says that if $R_{n_0}$ is sufficiently
concentrated around $1$ for some $n_0\ge0$, then it remains
concentrated for every $n>n_0$ and the free energy vanishes. In other
words, we establish a {\sl finite-volume condition for delocalization}.

\medskip

\begin{proposition}
  \label{th:fractmom} Let $B>2$ and $(\beta,h)$ be given. 
Assume that there exists $n_0\ge0$ and $(\log 2/\log B)<\gamma<1$ such that 
$\langle ([R_{n_0}-1]^+)^\gamma\rangle<B^\gamma-2$. Then, $\tf(\beta,h)=0$.
\end{proposition}

{\sl Proof of Proposition \ref{th:fractmom}.} We rewrite \eqref{eq:basic} as
\begin{equation}
\label{eq:basic2}
R_{n+1}-1 \, =\, \frac1B \left[\left( R_n^{(1)}-1 \right)
  \left( R_n^{(2)}-1 \right) +\left( R_n^{(1)}-1 \right)+
  \left( R_n^{(2)}-1 \right)\right],
\end{equation}
and we use the inequalities $[rs+r+s]^+ \le [r]^+[s]^+ +[r]^+ +
[s]^+$, that holds for $r, s \ge -1$, and $(a+b)^\gamma\le
a^\gamma+b^\gamma$, that holds for $\gamma \in (0,1]$ and $a, b \ge
0$.  If we set $A_n := \langle \left([ R_n -1]^+\right)^\gamma
\rangle$ we have
\begin{equation}
\label{eq:A}
A_{n+1}\, \le \, \frac1{B^\gamma} \left[A_n^2 + 2A_n\right]
\end{equation}
and therefore $A_n\searrow 0$ for $n\to\infty$ under the assumptions of the
Proposition.  Deducing $\tf(\beta,h)=0$ (and actually more than
that) is then immediate:
\begin{eqnarray}
  \langle \log R_n\rangle=\frac1\gamma  \langle \log (R_n)^\gamma\rangle
\le \frac1\gamma\langle \log \left[([R_n-1]^+)^\gamma+1\right]\rangle\le
\frac1\gamma\log (A_n+1)\stackrel{n\to\infty}\searrow 0.
\end{eqnarray}
\qed
\medskip

Proposition \ref{th:fractmom} will be essential in Section
\ref{sec:rel} to prove that, for $B>B_c$, an arbitrarily small
amount of disorder shifts the critical point. Let us also point out
that it implies that, if $\go_1$ is an unbounded random variable, then for any $B>2$ and
$\beta$ sufficiently large quenched and annealed critical points
differ (the analogous result for non-hierarchical pinning models was
proven in \cite[Corollary 3.2]{cf:T_fractmom}):
\medskip

\begin{cor} \label{cor:unbound}
Assume that $\bbP(\go_1>t)>0$ for every $t>0$. Then, for
every $h\in\R$ and $B>2$ there exists $\bar \beta_0<\infty$ such that
$\tf(\beta,h)=0$ for $\beta\ge\bar\beta_0$.
\end{cor}

\medskip
\noindent
{\it Proof of Corollary \ref{cor:unbound}}
Choose some $\gamma\in(\log 2/\log B,1)$.  One has $\lim_{\gb \to \infty}
  R_0=0$ $\bbP(\dd \go)$-a.s. (see \eqref{eq:R0} and note that
$\log M(\beta)/\beta\to\infty$ for $\beta\to+\infty$ under our assumption on
$\go_1$),
  while $\langle \left(([ R_0-1]^+)^\gamma \right)^{1/\gamma} \rangle
  \le 1+ \langle R_0 \rangle=1+\exp(h)$, so $\lim_{\gb \to \infty }A_0
  =0$.  \qed

 \medskip

 \begin{rem}\rm
 Note moreover that if 
 we set $X= \exp(\gb \go_1 -\log
 \M (\gb))$ we have (without requiring $\go_1$ unbounded) that $\langle
 ([(B-1) X -1]^+)^\gamma\rangle \stackrel{B\to \infty} \sim B^\gamma
 \langle X^\gamma \rangle$. The right-hand side is smaller than
 $B^\gamma -2$ for $X$ non-degenerate and $B$ large, so that if we
 choose $\gd>0$ such that $\exp( \gd \gamma) \langle X^\gamma \rangle
 <1$ we have
 \begin{equation}
 \langle ([(B-1)\exp(\gd) X -1]^+)^\gamma\rangle \, < \, B^\gamma-2, 
 \end{equation}
 for $B$ sufficiently large. Therefore, by applying
 Proposition~\ref{th:fractmom}, we see that for every $\gb>0$ there
 exists $\gd>0$ such that $\tf\left(\gb, h_c(0)+\delta\right)=0$ for
 $B$ sufficiently large. This observation actually follows also from
 the much more refined Proposition \ref{th:rel} below, which by the way
 says precisely how large $B$ has to be taken: $B>B_c$.
 \end{rem}  

 \begin{rem}\rm
\label{rem:as1}
  It follows from inequality \eqref{eq:A} that, if the assumptions of 
Proposition \ref{th:fractmom} are verified, then $A_n$ actually
vanishes exponentially fast for $n\to\infty$. Therefore, for $\gep>0$ one
has
\begin{equation}
\label{eq:BorCan}
  \bbP(R_n\ge 1+\gep)\, =\, \bbP([R_n-1]^+\ge\gep)\le \frac {A_n}{\gep^\gamma},
\end{equation}
and from the Borel-Cantelli lemma follows the almost sure convergence
of $R_n$ to $1$ when we recall that $R_n^{(i)}\ge r_n$ with $r_0=0$
($r_n$ is the solution of the iteration scheme \eqref{eq:r} and
converges to $1$).  
 \end{rem}

\subsection{Upper bounds on the free energy for $B >B_c$}

Here we want to prove the lower bound in \eqref{eq:eps_c}, plus the
fact that $h_c(\beta)>h_c$ whenever $\beta>0$ and $B>B_c$. This
follows from

\begin{proposition}
  \label{th:rel}
Let $B>B_c$. For every $\beta>0$ one has $h_c(\beta)>h_c(0)(= \log(B-1))$. 
Moreover,
there exists a positive constant $c$ (possibly depending on $B$)
such that for every $0\le \beta\le 1$
\begin{equation}
  \label{eq:rel}
  h_c(\beta) -h_c(0)
  \, \ge\, 
c \beta^{2\alpha/(2\alpha-1)}.
\end{equation}
\end{proposition}
Proposition \ref{th:rel} is proven in section \ref{sec:proofrel}, but
first we need to state a couple of technical facts.

\subsection{Auxiliary definitions and lemmas}
\label{sec:aux}
For $\gl \in \R$ and $N\in\N$ let $\bbP_{N, \gl}$ be defined by
\begin{equation}
\frac{
\dd \bbP_{N, \gl}}{\dd \bbP} (x_1,x_2, \ldots)\, =\,
\frac1 {\M(-\gl)^N} 
\exp\left(-\gl \sum_{i=1}^N x_i \right).
\end{equation} 

\smallskip

\begin{lemma}
\label{th:lemP}
There exists $1<C<\infty$ such that for $a\in(0,1)$, $\gd \in (0,
a/C)$ and $N\in\N$  we have
\begin{equation}
\label{eq:lemP}
\bbP_{N, \frac{\gd}{\sqrt{N}}}\left( \frac{\dd \bbP} {\dd \bbP_{N,
\frac{\gd}{\sqrt{N}}}} (\go)\, <\, \exp(-a) \right)\, \le \, C \left(
\frac{\gd}a\right)^2 .
\end{equation}
\end{lemma}
\smallskip

\noindent
{\it Proof of Lemma \ref{th:lemP}.}
We write
\begin{equation}
\label{eq:prtoest}
\bbP_{N, \frac{\gd}{\sqrt{N}}}\left( \frac{\dd \bbP} {\dd \bbP_{N,
\frac{\gd}{\sqrt{N}}}} (\go)\, <\, \exp(-a) \right)\, =\, \bbP_{N,
\frac\gd{\sqrt{N}}}\left( \gd \frac{\sum_{i=1}^{N} \go_i}{\sqrt{N}} +N \log
\M \left(-\frac{\gd }{\sqrt{N}}\right) \, < \, -a \right).
\end{equation}
Since all exponential moments of $\go_1$ are assumed to be finite, one has
\begin{equation}
 0\, \ge\,  \log \M(-\gl)-\gl\frac{\dd}{\dd \lambda}\left[\log \M(-\lambda)\right]\,\ge\,
 -\frac C2\gl^2,
\end{equation}
for some $1<C<\infty$ and $0\le\gl\le 1$ (the first inequality is due
 to convexity of $\gl\mapsto\log\M(-\gl)$). Note also that
\begin{equation}
\bbE_{N, \gl}(\go_1)\,=\,  
-\frac{\dd}{\dd \lambda}\left[\log \M(-\lambda)\right].
\end{equation}
Therefore, the right-hand side of \eqref{eq:prtoest}
is bounded above by
\begin{equation}
\bbP_{N, \gd/\sqrt{N}}\left(
 \frac{\sum_{i=1}^{N} \go_i}{\sqrt{N}} - \bbE_{N, \gd/\sqrt{N}}
\left[\frac{\sum_{i=1}^{N} \go_i}{\sqrt{N}}
\right]  <- \frac a{2\gd}
\right) \le \frac{4\delta^2}{a^2}\bbE_{N,\delta/\sqrt N}\left(
\go_1-\bbE_{N,\delta/\sqrt N}(\go_1)\right)^2,
\end{equation}
where we have used Chebyshev inequality and the fact that, under the
assumptions we made, $(a/\delta)-(C/2)\delta >a/(2\delta)$.  The proof
of \eqref{eq:lemP} is then concluded by observing that the variance of
$\go_1$ under $\bbP_{N,\gl}$ is $\dd^2/\dd \gl^2 \log \M(-\gl)$, which
is bounded uniformly for $0\le\gl\le 1$.  \qed
\bigskip

We define the sequence  $\{a_n\}_{n=0, 1, \ldots}$ by setting 
$a_0=a>0$ and $a_{n+1}= f(a_n)$ with
\begin{equation}
\label{eq:rec-inv}
f(x)\, :=\, \sqrt{B x+(B-1)^2} -(B-1). 
\end{equation}
We define also the sequence 
$\{b_n\}_{n=0,1, \ldots}$ by setting $b_0=b \in (-(B-2), 0)$ 
and $b_{n+1}= f(b_n)$.
Note that $a_n= g(a_{n+1})$ and $b_n= g(b_{n+1})$
for $g(x)=(2(B-1)x+x^2)/B$.

\smallskip

\begin{lemma} 
\label{th:ptlem}
There exist two constants $G_a>0$ et $H_b>0$ such that   for
$n \to \infty$
\begin{equation}
\label{eq:ptlem}
a_n\sim G_a \left(\frac{B}{2(B-1)}\right)^n\, =\, G_a 2^{-\ga n} 
\  \ \text{ and }\ \  b_n\sim -H_b \left(\frac{B}{2(B-1)}\right)^n
\, =\, -H_b 2^{-\ga n}
.
\end{equation}
Moreover, $G_a\stackrel{a\to0}\sim a$ and $H_b\stackrel{b\to0}\sim |b|$.
\end{lemma}
\smallskip

\noindent
{\it Proof of Lemma \ref{th:ptlem}.}
In order to lighten the proof we put $s:=B/(2(B-1))$ and we observe that
$0<s<1$ since $B>2$. The function
 $f(\cdot)$ is concave
  and  $f^\prime (0)=s$, so  $a_n$ vanishes exponentially fast:
\begin{equation}
 \label{eq:pop}
a_n\le a\, s^n.
\end{equation}
  Moreover, 
\begin{equation}
\label{eq:monoton}
  \frac{a_n}{s^n}=\frac{a_{n-1}}{s^{n-1}}\frac1{1+a_n/(2(B-1))}\, \ge\,
\frac{a_{n-1}}{s^{n-1}}\frac1{1+a s^n/(2(B-1))},
\end{equation}
so that for every $n>0$
\begin{eqnarray}
\label{eq:Ca2}
\frac{a_n}{s^n}\ge a\, \prod_{\ell=1}^\infty\frac1{1+a s^\ell/(2(B-1))}>0.
\end{eqnarray}
From \eqref{eq:monoton} we see that $a_n\,s^{-n}$ is monotone increasing in 
$n$, so that the first statement in \eqref{eq:ptlem} holds with
$G_a\in (0,a)$ by  \eqref{eq:pop} and \eqref{eq:Ca2}.
The fact that $G_a\sim a$ for $a\to0$ follows from the fact that
the product in  \eqref{eq:Ca2} converges to $1$ in this limit.
\smallskip

The second relation is proven in a similar way. Since $b_n<0$ for every
$n$, one has first of all
\begin{eqnarray}
\label{eq:Cb}
  \frac{b_n}{s^n}=\frac{b_{n-1}}{s^{n-1}}\frac1{1+b_n/(2(B-1))}<
\frac{b_{n-1}}{s^{n-1}}.
\end{eqnarray}
Moreover, since $\vert b_n \vert $ decreases to zero and $f(x)\ge
c_1(b) x$ for $b\le x\le0$ for some $c_1(b)<1$ if $b>-(B-2)$, one
sees that $|b_n|$ actually vanishes exponentially fast.  Therefore,
from \eqref{eq:Cb}
\begin{eqnarray}
    \frac{b_n}{s^n}\ge \frac{b_{n-1}}{s^{n-1}}\frac1{1-c_2(b)\,c_1(b)^n}
\ge b\prod_{\ell=1}^\infty\frac1{1-c_2(b)\,c_1(b)^\ell}.
\end{eqnarray} 
One has then the second statement of \eqref{eq:ptlem} with
$H_b\in(|b|,\infty)$.

\qed

\subsection{Proof of Proposition \ref{th:rel}}
\label{sec:proofrel}
In this proof $C_i,i=1,2,\ldots$ denote constants depending only on $\gb_0$ and (possibly)
on $B$.
Recall that the exponent $\alpha$ defined in \eqref{eq:alpha}
satisfies $1/2<\ga<1$ for $B>B_c$. Fix $\beta_0>0$, let
$0<\beta<\beta_0$ and choose $h=h(\beta)$ such that
\begin{eqnarray}
  \med{R_0}=(B-1)+\eta\gb^{\frac{2\ga}{2\ga-1}},
\end{eqnarray}
where $\eta>0$ will be chosen sufficiently small and independent of
$\beta$ later.  Call $n_0:=n_0(\eta,\gb)$ the integer such that
\begin{eqnarray}
\label{eq:prima}
  \med{R_{n_0}}\le B\le \med{R_{n_0+1}},
\end{eqnarray}
{\sl i.e.}, $P_{n_0}\le 1\le P_{n_0+1}$.  
Note that
$n_0(\eta, \gb)$ becomes larger and larger
as $\gb\searrow 0$: this can be quantified since 
from  \eqref{eq:P} one sees that
$a_n:=P_{n_0-n}$ satisfies for $0\le n< n_0$ the iteration
$a_{n+1}=f(a_n)$ introduced in \S~\ref{sec:aux},
 and therefore it
follows from Lemma \ref{th:ptlem} that 
\begin{equation}
\label{eq:n0}
  \left|n_0(\eta,\gb)-\log\left(\eta^{-1}\gb^{-\frac{2\ga}{2\ga-1}}\right)
/(\alpha\log 2)\right|\, \le\,  C_1,
\end{equation}
for every $0<\eta<1/C_1$ and $\beta\in[0,\beta_0]$.  With the 
notations of Section  \ref{sec:aux}, let $\tilde
\bbP:=\bbP_{2^{n_0},\delta 2^{-n_0/2}}$, where $\delta:=\delta(\eta)$ will
be chosen suitably small later.
Note that, with $\gl:=\delta 2^{-n_0/2}$, one has from \eqref{eq:n0}
\begin{eqnarray}
\label{eq:lpetit}
  \frac1{C_2}\delta\,\eta^{1/(2\alpha)}\gb^{1/(2\ga-1)} \le \gl\le
C_2\delta\,\eta^{1/(2\alpha)}  \gb^{1/(2\ga-1)}.
\end{eqnarray}
In particular, since $\alpha<1$,  if $\eta$ is small enough
then $\gl\le \beta$ uniformly for $\beta\le \beta_0$. Observe also that
\begin{eqnarray}
\label{eq:R0'}
  \tilde \bbE(R_0)=\med{R_0}\frac{\M(\gb-\gl)}{\M(\gb)\M(-\gl)},
\end{eqnarray}
and call $\phi(\cdot):=\log \M(\cdot)$. Since $\phi(\cdot)$ is strictly convex,
one has
\begin{eqnarray}
  \phi(\gb-\gl)-\phi(\gb)-\phi(-\gl)=-\int_{-\gl}^0 \dd x\int_0^\gb\dd y
\,\phi''(x+y)\in \left(- \frac{\gl\gb}{C_3},-C_3\gl\gb\right),
\end{eqnarray}
for some $C_3>0$, uniformly in $\beta\le \beta_0$ and $0\le \lambda\le
\beta$  and, thanks to \eqref{eq:lpetit}, if $\eta$ is chosen sufficiently
small,
\begin{eqnarray}
1-\frac{\beta\gl}{C_4} \le \frac{\M(\gb-\gl)}{\M(\gb)\M(-\gl)}\le 1- C_4\beta\gl.
\end{eqnarray}
Therefore, from  \eqref{eq:R0'} and \eqref{eq:lpetit} and choosing
\begin{eqnarray}
  \label{eq:deltall}
\eta^{1-1/(2\ga)}\ll \delta(\eta)\ll 1,
\end{eqnarray}
(which is possible with $\eta$ small since $\alpha>1/2$) one has
\begin{eqnarray}
-C_5^{-1}\,\delta(\eta)\,
    \eta^{1/(2\alpha)}\beta^{\frac{2\alpha}{2\alpha-1}}< \tilde
    \bbE(R_0)-(B-1)\le-C_5\delta(\eta)\,
    \eta^{1/(2\alpha)}\beta^{\frac{2\alpha}{2\alpha-1}},
\end{eqnarray}
always uniformly in $\beta\le \beta_0$.

Since $b_n:=\tilde\bbE (R_{n_0-n})-(B-1)$ satisfies the recursion
$b_{n+1}=f(b_n)$, from the second statement of  \eqref{eq:ptlem} if
follows that
\begin{equation}
  \tilde \bbE R_{n_1}\, \le\,  \frac {B}2,
\end{equation}
for some integer $n_1:=n_1(\eta,\beta)$ satisfying
\begin{equation}
  \label{eq:n1} 
n_1\, \le\, \left(
\log\left(
\delta(\eta)^{-1}\eta^{-1/(2\alpha)}\beta^{-2\ga/(2\ga-1)}
\right)
\Big /
\left(
\ga \log 2\right) \right)
+C_6.
\end{equation}
It is immediate to see that $n_0(\eta,\beta)-n_1(\eta,\beta)$ gets
 large (uniformly in $\beta$) for $\eta$ small, if condition
\eqref{eq:deltall} is satisfied. Therefore, since the fixed point $1$
of the iteration for $\tilde \bbE R_{n}$
is attractive, one has that
\begin{equation}
   \tilde \bbE R_{n_0}\le 1+r_1(\eta),
\end{equation}
(here and in the following, $r_i(\eta)$ with $i\in\N$ denotes a positive
quantity which vanishes for $\eta\searrow0$, uniformly in $\gb\le \gb_0$).
On the other hand, 
one has deterministically 
\begin{eqnarray}
\lim_{n\to\infty}  [1-R_n]^+=0,
\end{eqnarray}
as one sees immediately comparing the evolution of $R_n$ with that
obtained setting $R_0^{(i)}=0$ for every $i$. In particular,
$R_{n_0}\ge 1-r_2(\eta)$.  An
application of Markov's inequality gives
\begin{eqnarray}
  \tilde \bbP(R_{n_0}\ge 1+r_3(\eta))\le r_3(\eta).
\end{eqnarray}
It is immediate to prove that, given a random variable $X$ and two 
mutually absolutely continuous laws $\bbP$ and $\tilde \bbP$, one has
for every $x,y>0$
\begin{eqnarray}
  \bbP(X\le 1+x)\ge e^{-y}\left[\tilde\bbP(X\le 1+x)-\tilde \bbP\left(
\frac{\dd\bbP}{\dd\tilde\bbP}\le e^{-y}\right)\right].
\end{eqnarray}
Applying this to the case $X=R_{n_0}$ and using Lemma \ref{th:lemP} with 
$r_4(\eta)>C\delta(\eta)$ gives
\begin{eqnarray}
  \bbP(R_{n_0}\le 1+r_3(\eta))
\ge e^{-r_4(\eta)}\left[1-r_3(\eta)-C\left(\frac{\delta(\eta)}{r_4(\eta)}
\right)^2\right].
\end{eqnarray}
In particular,
choosing 
\begin{equation}
\delta(\eta)\ll r_4(\eta)\ll1,
\end{equation}
one has 
\begin{eqnarray}
\label{eq:seconda}
  \bbP(R_{n_0}\le 1+r_3(\eta))\ge 1-r_5(\eta),
\end{eqnarray}
and we emphasize that this inequality holds uniformly in $\gb\le
\gb_0$.

At this point  \eqref{eq:rel} is essentially proven: choose some
$\gamma\in(\log 2/\log B,1)$ and observe that
\begin{equation}
\begin{split}
 \med{\left([R_{n_0}-1]^+\right)^\gamma} \,&\le\, 
r_3(\eta)^\gamma+\left(\bbE[R_{n_0}-1]^+\right)^\gamma
\left(\bbP(R_{n_0}\ge 1+r_3(\eta))\right)^{1-\gamma}\\
\label{eq:uff}
&\le\, r_3(\eta)^\gamma+B^\gamma r_5(\eta)^{1-\gamma},
\end{split}
\end{equation}
where in the first inequality we have used H\"older inequality and in the
second one we have used  \eqref{eq:prima} and \eqref{eq:seconda}.  Finally, we
remark that the quantity in \eqref{eq:uff} can be made smaller than
$B^\gamma-2$ choosing $\eta$ small enough.  At this point, we can
apply Proposition \ref{th:fractmom} to deduce that $\tf(\beta,h)=0$
for $h=\log (B-1)+\eta \beta^{2\ga/(2\ga-1)}$ with $\eta$ small but
finite, which proves  \eqref{eq:rel}.

We complete the proof by observing 
 that $h_c(\beta)>\log (B-1)$ for every $\beta>0$ follows from
the arbitrariness of $\beta_0$.  \qed

\section{The delocalized phase}
\label{sec:deloc}
Here we prove Theorem \ref{th:paths} using the
representation \eqref{eq:EBgen}, given in Appendix \ref{sec:stand}, for $R_n$.
With reference to \eqref{eq:EBgen}, let us observe that
\begin{equation}
\label{eq:fewcontacts}
\lim_{n\to\infty} p(n,\emptyset)\,=\,1,
 \end{equation}
which is just a way to interpret 
\begin{equation}
 \lim_{n\to\infty}r_n\,=\,1.
\end{equation}
when $r_0=0$, that follows directly from \eqref{eq:r}.



\medskip


Fix $\epsilon>0$ arbitrarily small and consider $h<h_c(\gb)$.  Let $\bar R_n$ be the
partition function which corresponds to $h_c(\gb)$ and $R_n$ the one
that corresponds to $h$.  We can find $K$ large enough such that
\begin{equation}\label{eq:barR}
 \bbP \left(\bar R_n \ge K \right)\le \epsilon/2 \quad \text{ for all } n\ge 1.
\end{equation}
This follows from the fact that 
$\bar R_n\ge (B-1)/B$, and from \eqref{eq:nondecrease}.
We define $C:=(\log(2K/\epsilon))/{(h_c(\beta)-h)}$ and we write, using \eqref{eq:EBgen},
\begin{multline}
  R_n\, =\,p(n,\emptyset)+
  \sumtwo{\mathcal I \subset\{1,\dots,2^n\}}{1\le |\mathcal I|\le C}
p(n,\mathcal I) \exp\left(\sum_{i\in \mathcal I}(\gb\go_i-\log M(\gb)+h)\right)
\\
  +\sumtwo{\mathcal I \subset\{1,\dots,2^n\}}{|\mathcal I|
    >C}p(n,\mathcal I) \exp\left(\sum_{i\in \mathcal I}(\gb\go_i-\log
    M(\gb)+h)\right)\, =:\, T_1+T_2+T_3.
\end{multline}
$T_1$ is smaller than $1$ and 
\begin{equation}
T_3\, \le \, 
 \exp\left(-C(h_c(\beta)-h)\right)\bar R_n,
\end{equation}
so that $T_3 \le \epsilon/2$ with probability greater than $(1-\epsilon/2)$
({\sl cf.} \eqref{eq:barR}) for all $n$.  As for $T_2$, its easy to
compute and bound its expectation:
\begin{equation}
  \left\langle \sumtwo{\mathcal I \subset\{1,\dots,2^n\}}
    {1\le |\mathcal I|\le C}p(n,\mathcal I) \exp\left(\sum_{i\in \mathcal I}(\gb\go_i-\log M(\gb)+h)\right)\right\rangle \, \le\,  \exp(Ch)[1-p(n,\emptyset)],
\end{equation}
and \eqref{eq:fewcontacts} tells us that the right-hand side tends to
zero when $n$ goes to infinity.  In particular we can find $N$
(depending on $C$) such that for all $n\ge N$ we have
\begin{equation}
  \left\langle\sumtwo{\mathcal I \subset\{1,\dots,2^n\}}
{1\le |\mathcal I|\le C}p(n,\mathcal I) \exp\left(\sum_{i\in \mathcal I}(\gb\go_i-\log M(\gb)+h)\right)\right\rangle\, \le\,  \epsilon^2/4.
\end{equation}
Then for $n\ge N$ we have
$
 \bbP(T_2\ge \epsilon/2)\le \epsilon/2$.
Altogether we have
\begin{equation}
 \bbP(R_n\ge 1+\epsilon)\,\le\,  \epsilon,
\end{equation}
and since $R_n$ is bounded from below by $p(n,\emptyset)$ which tends to $1$,
the proof is complete.
\qed

\appendix

\section{Existence of the free energy and annealed system estimates}

\label{sec:stand}

\subsection{Proof of Theorem \ref{th:F}}


Since the basic induction \eqref{eq:basic} gives $ R_n\ge (B-1)/B $
for every $n\ge 1$, one has
\begin{equation}
\frac{R_{n+1}}{B}\, \ge\,  \frac{R_n^{(1)}}{B}\frac{R_n^{(2)}}{B}, \label{eq:ge}
\end{equation}
and
\begin{equation}
R_{n+1} \, \le\,  \frac{R_n^{(1)}R_n^{(2)}}{B}+ \frac{B}{B-1}R_n^{(1)}R_n^{(2)},
\end{equation}
so  that
\begin{align}
  (K_B R_{n+1}) \le (K_BR_n^{(1)})(K_BR_n^{(2)}) \quad \text{with}
  \quad K_B=\frac{B^2+B-1}{B(B-1)} \label{eq:le}.
\end{align}
Taking the logarithm of \eqref{eq:ge} and \eqref{eq:le}, we get that
\begin{equation}
\label{eq:nondecrease}
\left\{2^{-n}\E\big[\log(R_n/B)\big]\right\}_{n=1, 2, \ldots}
\ \text{
 is non-decreasing,} 
 \end{equation}
 while
 \begin{equation}
\label{eq:nonincrease}
\left\{2^{-n}\E\big[\log(K_B R_n)\big]\right\}_{n=1, 2, \ldots} \
\text{ is non-increasing},
\end{equation}
so that both sequences are converging to the same limit
\begin{equation}
\tf(\beta,h)=\lim_{n\rightarrow\infty}2^{-n}\med{\log R_n}
\end{equation}
and \eqref{eq:encadre} immediately follows.  It remains to be proven that
the limit of $2^{-n}\log R_n$ exists $\bbP(\dd \go)$--almost surely and  in $\bbL^1(\dd \bbP)$.
Fixing some $k\ge1$ and iterating \eqref{eq:ge} one obtains for $n>k$
\begin{equation}
  2^{-n}\log (R_{n}/B)\ge 2^{-k}\big(2^{k-n}\sum_{i=1}^{2^{n-k}}\log
  (R_k^{(i)}/B)\big).
\end{equation}
Using the strong law of large numbers in the right-hand side, we get
\begin{equation}
  \liminf_{n\to\infty} 2^{-n} \log (R_n/B) \ge 2^{-k}\med{\log
  (R_k/B)} \quad \bbP(\dd \go)- \text{a.s.}.
\end{equation}
Hence taking the limit for $k\rightarrow\infty$ in the right-hand side
again we obtain
\begin{equation}
  \liminf_{n\to\infty} 2^{-n} \log R_n=\liminf_{n\to\infty} 2^{-n}
  \log (R_n/B) \ge \tf(\beta,h) \quad \bbP(\dd \go)- \text{a.s.}.
\end{equation}
Doing the same computations with \eqref{eq:le} we obtain
\begin{equation}
  \limsup_{n\to\infty} 2^{-n} \log R_n=\limsup_{n\to\infty} 2^{-n} \log (K_B R_n)
  \le \tf(\beta,h) \quad \bbP(\dd \go)- \text{a.s.}.
\end{equation}
This ends the proof for the almost sure convergence. The proof of the
$\bbL^1(\dd \bbP)$ convergence is also fairly standard, and we leave it to
the reader.

The fact that $\tf(\beta,\cdot)$ is non-decreasing follows from the
fact that the same holds for $R_n(\beta,\cdot)$, and this is easily
proved by induction on $n$. Convexity of $(\gb, h) \mapsto \tf (\gb, h
+ \log \M(\gb))$ is immediate from \eqref{eq:DHV} (hence for $B=2,3,
\ldots$). But \eqref{eq:DHV} can be easily generalized to every $B>1$:
this follows by observing that from \eqref{eq:basic} and \eqref{eq:R0}
one has that
\begin{equation}
\label{eq:EBgen}
 R_n\, =\, \sum_{\mathcal{I} \subset\{1,\dots,2^n\}}p(n,\mathcal I) \exp\left(\sum_{i\in \mathcal I}(\gb\go_i-\log M(\gb)+h)\right),
\end{equation}
for suitable positive values $p(n, \mathcal I)$, which depend on $B$:
by setting $\gb=h=0$ we see that $ \sum_{\mathcal{I}}p(n,\mathcal
I)=1$ and hence $R_n$ can be cast in the form of the expectation of a
Boltzmann factor, like \eqref{eq:DHV}. This yields the desired
convexity.  \qed

\smallskip

\begin{rem}
\label{rem:5} 
\rm
Another consequence of \eqref{eq:EBgen} is that
 $\tf(\gb, h+ \log \M( \gb))\ge \tf (0,  h)$ \cite[Ch.~5, Prop.~5.1]{cf:Book}.
\end{rem}


\medskip

\subsection{Proof of Theorem \ref{th:pure}} 
When $\gb=0$  the iteration \eqref{eq:basic} reads
\begin{equation}
R_{n+1}=\frac{R_n^2+(B-1)}{B}.
\end{equation}
A quick study of the function $x\mapsto [x^2+(B-1)]/B$, gives that
$R_n\stackrel{n\rightarrow\infty}{\rightarrow}\infty$ if and only if
$R_0>(B-1)$. Initial conditions $R_0<B-1$ are attracted by the stable
fixed point $1$, while the fixed point $(B-1)$ is unstable.  The
inequality \eqref{eq:ge} guaranties that $\tf(0,h)>0$ when $R_N>B$ for
some $N$. This immediately shows that that $h_c(0)=\log(B-1)$.


Next we prove \eqref{eq:alpha0}, {\sl i.e.}, that (with the notations in
\eqref{eq:gep} and \eqref{eq:hatF}) there exists a constant $C$
such that 
\begin{equation}
\label{eq:modif}
\frac{1}{C}\gep^{1/\ga}\, \le\,  \hat{\tf}(0,\gep) \, \le\,  C \gep^{1/\ga} 
\end{equation}
for all $\gep\in(0,1)$.  To that purpose take $a:=a_0$ such that
$\hat{\tf}(0,a)=1$ 
(this is possible because of the  convexity  of
$\tf(\gb, \cdot+ \log \M(\gb))$ we obtain both continuity and 
$\lim_{a \to \infty}\hat{\tf}(0,a)=\infty$) 
and note that  the sequence $\{a_n\}_{n\ge0}$ defined
just before Lemma \ref{th:ptlem} 
is such that $2\,\hat{\tf}(0,a_{n+1})=\hat{\tf}(0,a_{n})$, so that
$\hat \tf(0,a_{n+1})=2^{-n}$.
Thanks to Lemma \ref{th:ptlem} we have that along this sequence
\begin{equation}
\hat{\tf}(0,a_n)\sim  2G_a^{-1/\ga} a_n^{1/\ga}. \label{eq:sim}
\end{equation}
Let $K_a$ be such that $a_n\le K_a a_{n+1}$ for all $n$, and $c_a$
such that $c_a^{-1}a_n^{1/\ga}\le\hat{\tf}(0,a_n)\le
c_a\,a_n^{1/\ga}$.  Then, for all $n$ and all
$\gep\in[a_{n+1},a_{n}]$, since $\hat{\tf}(0,\cdot)$ is increasing we
have
\begin{equation}
\begin{split}
  \hat{\tf}(0,\gep)&\ge \hat{\tf}(0,a_{n+1}) 
\ge c_a^{-1}a_{n+1}^{1/\ga}\ge c_a^{-1}K_a^{-1/\ga}\gep^{1/\ga},\\
  \hat{\tf}(0,\gep)&\le \hat{\tf}(0,a_{n})\le c_a a_n^{1/\ga}\le
  c_aK_a^{1/\ga}\gep^{1/\ga}.
  \end{split}
\end{equation}
Finally, the analyticity of $\tf(0, \cdot)$ on $(h_c, \infty)$  follows 
for example from \cite[Lemma~4.1]{cf:OvKr}.

\qed

\subsection{About   models with $B\le 2$}
\label{sec:Bsmallerthan2}

We have chosen to work with the model \eqref{eq:R}, with positive
initial data and $B>2$, because this is the case that is directly
related to pinning models and because in this framework we had the
precise aim of proving the physical conjectures formulated in
\cite{cf:DHV}. But of course the model is well defined for all $B\neq
0$ and in view of the direct link with the logistic map  $z \mapsto A
z(1-z)$, {\sl cf.}  \eqref{eq:logistic}, also the case $B\le 2$ appears to be
intriguing. Recall that $A=2(B-1)/B$ and note that $A \in (1,2)$ if
$B\in (2, \infty)$.  What we want to point out here is mainly that the
case of \eqref{eq:R} with positive initial data and $B \in (1,2)$,
{\sl i.e.} $A\in (0,1)$, is already contained in our analysis.  This
is simply the fact that there is a duality transformation relating
this new framework to the one we have considered. 
Namely, if we let $B\in(1,2)$ and we set $\hat R_n := R_n/(B-1)$, then
$\hat R_n$ satisfies \eqref{eq:R} with $B$ replaced by $\hat B:=
B/(B-1) >2$. Of course the fixed points of $x\mapsto (x^2 +\hat B -1)/
\hat B$ are again $1$ (stable) and $\hat B -1$ (unstable). This
transformation allows us to generalize immediately all the theorems we
have proven in the obvious way, in particular the marginal case
corresponds to $\hat B=\hat B_c:=2+\sqrt{2}$, {\sl i.e.}, $B= \sqrt{2}$
and in the irrelevant case ($B \in(\sqrt{2},2)$) the condition on $\gD (\gb)$
in Theorem~\ref{th:eps_c} now 
reads $ \M(2\gb)/(\M(\gb))^2 < B^2 -1$ 
  
 This discussion leaves open the cases $B=1$ and $B=2$ to
 which we cannot apply directly our theorems, but:
 \begin{enumerate}
 \item If $ B=1$ the model is exactly solvable and $R_n$ is equal to
   the product of $2^n$ positive IID random variables distributed like
   $R_0$, so $\tf (\gb, h)= h-\log \M (\gb)$. The model in this case
   is a bit anomalous, since the stable fixed point is $0$ and
   therefore the free energy can be negative and no phase transition
   is present (this appears to be the analogue of the non-hierarchical
   case with inter-arrival probabilities that decay exponentially fast
   \cite[Ch.~1, Sec.~9]{cf:Book}).
\item If $B=2$ then, with reference to \eqref{eq:r}, $r_n \nearrow
  \infty$ if $r_0>1$ and $r_n \nearrow 1$ if $r_0<1$. The basic
  results like Theorem~\ref{th:F} are quickly generalized to cover this
  case. Only slightly more involved is the generalization of the other
  results, notably Theorem~\ref{th:eps_c}$(1)$.  In fact we cannot
  apply directly our results because the iteration for $P_n$, that is
  $\Rav{n} -1$, reads $P_{n+1}=P_n +(P_n^2 /2)$ ({\sl cf.}
  \eqref{eq:P}) so that the growth of $P_n$, for $P_0>0$, is just due
  to the nonlinear term and it is therefore slow as long as $P_n$ is
  small.  However the technique still applies (note in particular that, by 
   \eqref{eq:Q} and \eqref{eq:QboundP}, the variance of $R_n$
  decreases exponentially if $\gD_0<2$ as long as $P_n$ is sufficiently
  small)  and along this line one shows that the disorder is
  irrelevant, at least as long as $\gD_0 <2$.
 \end{enumerate} 

\medskip

If we now let $B $ run from $1$ to infinity, we simply conclude that
the disorder is irrelevant if $B \in (\sqrt{2}, 2+\sqrt{2})$, 
and it is instead relevant in $B \in (1,
\sqrt{2}) \cup (2+\sqrt{2}, \infty)$.  In the case $B=1$ (and, by
duality, $B=\infty$) there is no phase transition.

\medskip

Finally, a word about the models with $B<1$.  Various cases should be
distinguished: going back to the logistic map, we easily see that
playing on the values of $B$ one can obtain values of $\vert A\vert $
larger than $2$ and the very rich behavior of the logistic map sets
in \cite{cf:Beardon}: non-monotone convergence to the fixed point, oscillations in a
finite set of points, chaotic behavior, unbounded trajectories for any
initial value.  It appears that it is still  possible to generalize
our approach to deal with some of these cases, but this would lead us
far from our original aim.
Moreover, for $B<1$ the property of positivity of $R_n$, and
therefore its statistical mechanics interpretation as a partition
function, is lost.

\section*{acknowledgments} We are very much indebted to Bernard Derrida
for drawing our attention to the hierarchical pinning model and for
several enlightening discussions. We acknowledge also the support of
ANR grant {\sl POLINTBIO}, and F.L.T. acknowledges the support of ANR grant
{\sl LHMSHE}.

\bigskip

\noindent
{\bf Note added in proof.}
After this work was completed, a number of results have been proven by
developing further the ideas set forth here, solving some of questions
raised at the end of Section~\ref{sec:pinning}. First of all we were
able (in collaboration with B. Derrida \cite{cf:DGLT}) to extend the
main idea of this work to the non hierarchical set-up and we have
shown that the quenched critical point (of the non-hierarchical model)
is shifted with respect to the annealed value for arbitrarily small
disorder, if $\ga>1/2$ (this result has been sharpened in
\cite{cf:AZ}, taking a different approach).  Then one of us
\cite{cf:Hubert} has been able to show the shift of the critical point
for arbitrarily small disorder for $\ga =1/2$ in a hierarchical with
site disorder (the case considered here is  bond disorder, {\sl
  cf.}  Figure~\ref{fig:diamond}) by using a location-dependent shift of
the disorder variables in the change-of-measure argument (in the present paper,
the shift is the same for each variable). Finally, very recently
\cite{cf:GLT_marg} we have also been able to treat the case $\ga =1/2$
($B=B_c$), both for the hierarchical and non-hierarchical model, by
introducing long range correlations in the auxiliary measure 
$\tilde \bbP$.

\end{document}